\documentclass[reqno,oneside]{amsart} 

\usepackage{amsthm,
amssymb,amsmath,enumerate,stackrel,
verbatim,
}
\usepackage{xr-hyper}
\usepackage[all,cmtip]{xy}
\usepackage{enumitem}

\usepackage[pagebackref,hyperindex,linktocpage=true]{hyperref}
\hypersetup{
    colorlinks,
    linkcolor={red!50!black},
    citecolor={red!50!black}, 
    urlcolor={red!50!black}, 
    filecolor={red!50!black} 
}

\usepackage{tikz}
\usetikzlibrary{matrix,arrows}

\usepackage{tikz-cd}

\def\label#1{\label{#1}}

\usepackage{theoremref}
\definecolor{labelkey}{rgb}{1,0,0}

\makeatletter
\@addtoreset{equation}{section}
\makeatother

\numberwithin{equation}{section}

\theoremstyle{definition}
\newtheorem{Defi}{Definition}[section] \newcommand{\defi}{\begin{Defi}} \newcommand{\xdefi}{\end{Defi}} 
\newtheorem{DefiLemm}[Defi]{Definition and Lemma} \newcommand{\defilemm}{\begin{DefiLemm}} \newcommand{\xdefilemm}{\end{DefiLemm}} 
\newtheorem{Bsp}[Defi]{Example} \newcommand{\exam}{\begin{Bsp}} \newcommand{\xexam}{\end{Bsp}} 
\newtheorem{Syno}[Defi]{Synopsis} \newcommand{\syno}{\begin{Syno}} \newcommand{\xsyno}{\end{Syno}} 
\newtheorem{Bem}[Defi]{Remark} \newcommand{\rema}{\begin{Bem}} \newcommand{\xrema}{\end{Bem}} 
\newtheorem{Notation}[Defi]{Notation} \newcommand{\nota}{\begin{Notation}} \newcommand{\xnota}{\end{Notation}} 
\newtheorem{Convention}[Defi]{Convention} \newcommand{\conv}{\begin{Convention}} \newcommand{\xconv}{\end{Convention}} 
\newtheorem{Warning}[Defi]{Warning} \newcommand{\warn}{\begin{Warning}} \newcommand{\xwarn}{\end{Warning}} 
\newtheorem{Situation}[Defi]{Situation} \newcommand{\situ}{\begin{Situation}} \newcommand{\xsitu}{\end{Situation}} 

\theoremstyle{plain}
\newtheorem{Theo}[Defi]{Theorem} \newcommand{\theo}{\begin{Theo}} \newcommand{\xtheo}{\end{Theo}} 
\newtheorem{Satz}[Defi]{Proposition} \newcommand{\prop}{\begin{Satz}} \newcommand{\xprop}{\end{Satz}} 
\newtheorem{Lemm}[Defi]{Lemma} \newcommand{\lemm}{\begin{Lemm}} \newcommand{\xlemm}{\end{Lemm}} 
\newtheorem{Coro}[Defi]{Corollary} \newcommand{\coro}{\begin{Coro}} \newcommand{\xcoro}{\end{Coro}}
\newtheorem{Ques}[Defi]{Question} \newcommand{\ques}{\begin{Ques}} \newcommand{\xques}{\end{Ques}}
\newtheorem{Conj}[Defi]{Conjecture} \newcommand{\conj}{\begin{Conj}} \newcommand{\xconj}{\end{Conj}}

\newcommand{\refsect}[1]{Section \ref{sect--#1}}
\newcommand{\refit}[1]{(\ref{item--#1})}
\newcommand{\refeq}[1]{(\ref{eqn--#1})}

\newcommand{\eqn}{\begin{equation}} \newcommand{\xeqn}{\end{equation}}
\newcommand{\eqnarr}{\begin{eqnarray*}} \newcommand{\xeqnarr}{\end{eqnarray*}}
\newcommand{\eqnarra}{\begin{eqnarray}} \newcommand{\xeqnarra}{\end{eqnarray}}

\newcommand{\pf}{\begin{proof}} \newcommand{\xpf}{\end{proof}}

\newcounter{heyheyCounter}[section]

\newcommand{\nc}{\newcommand}
\nc{\StP}[1]{\cite[Tag~\href{http://stacks.math.columbia.edu/tag/#1}{#1}]{StacksProject}}
\nc{\StPd}[2]{\cite[Tags~\href{http://stacks.math.columbia.edu/tag/#1}{#1}, \href{http://stacks.math.columbia.edu/tag/#2}{#2}]{StacksProject}} 

\nc{\on}{\operatorname}
\nc{\aff}{{\on{aff}}}
\nc{\modi}{{\on{mod}}} 
\nc{\even}{{\on{even}}}
\nc{\odd}{{\on{odd}}}
\nc{\naive}{{\on{naive}}}
\nc{\hofib}{\on{hofib}}
\nc{\Bun}{\on{Bun}}
\nc{\ad}{{\on{ad}}}
\nc{\lft}{{\on{lft}}}

\nc{\Weil}{{\on{Weil}}} 
\nc{\FWeil}{{\on{FWeil}}} 
\nc{\cons}{{\on{cons}}} 
\nc{\tot}{{\on{Tot}}} 

\nc{\str}{\on{-}}
\nc{\perf}{{\on{perf}}}
\nc{\Rel}{{\on{Pos}}}
\nc{\lan}{\langle}
\nc{\ran}{\rangle}

\nc{\PM}{{\on{PM}}}
\nc{\calH}{{\mathcal H}}
\nc{\calO}{{\mathcal O}}
\nc{\calM}{{\mathcal M}}
\nc{\co}{\colon}

\newcommand{\category}[1]{\mathrm{#1}}

\def\gr{\category {gr}} 
\newcommand{\Ho}{\category{Ho}} 
\newcommand{\Mod}{\category{Mod}} 
\newcommand{\Modf}{\Mod^{\operatorname{f}}} 
\newcommand{\Shv}{\category{Shv}} 
\newcommand{\Perf}{\category{Perf}} 
\newcommand{\grPerf}{\category{grPerf}} 

\newcommand{\ModL}{\Mod_\Lambda} 
\newcommand{\SgrMod}{\category{grSMod}} 

\renewcommand{\Pr}{\category{Pr}}
\newcommand{\PrSt}{\Pr^{\category{St}}} 
\newcommand{\PrStL}{\PrSt_\Lambda} 
\newcommand{\grPrStL}{\gr\PrStL} 
\newcommand{\Sp}{\category{Sp}} 
\newcommand{\Cat}{\category{Cat}} 
\newcommand{\Fun}{\category{Fun}} 
\newcommand{\Corr}{\category{Corr}} 
\newcommand{\Sch}{\category{Sch}} 
\newcommand{\SmAff}{\category{SmAff}} 
\newcommand{\IndSch}{\category{IndSch}} 
\newcommand{\Alg}{\category{Alg}} 
\def\ft{\mathrm{ft}} 
\newcommand{\Fin}{\category{Fin}} 
\newcommand{\dual}{\vee} 
\newcommand{\Zar}{\mathrm{Zar}} 

\def\proper{\mathrm {proper}} 
\def\sep{\mathrm {sep}} 
\def\all{\mathrm {all}} 
\def\W{\mathrm {W}} 
\def\wt{\mathrm {w}} 
\def\tt{\mathrm {t}} 
\newcommand{\grMod}{\gr\Mod} 
\newcommand{\grModL}{\grMod_\Lambda} 
\def\Gm{\mathbf {G}_\mathrm m} 
\newcommand{\GmX}[  1]{\mathbf {G}_{\mathrm {m}, #1}} 

\newcommand{\colim}{\operatornamewithlimits{colim}} 

\def\id{{\rm id}} 
\def\ev{{\operatorname {ev}}} 
\def\pr{{\rm pr}} 
\def\opp{{\rm op}} 
\def\To#1#2{\mathop{\count0=#1 \loop\ifnum\count0>0 \smash-\mkern-7mu \advance\count0 -1 \repeat \mathord\rightarrow}\limits^{#2}} 
\def\Maps{\mathop{\rm Maps}\nolimits} 
\def\IMaps{\underline {\Maps}} 
\def\Char{\mathop{\rm char}\nolimits} 
\def\CH{\mathop{\rm CH}\nolimits} 
\def\laxlim{\mathop{\rm laxlim}} 
\def\Hom{\mathop{\rm Hom}\nolimits} 
\def\fib{\mathop{\rm fib}\nolimits} 

\def\Ind{\mathop{\category{Ind}}} 
\def\Gr{\mathop{\rm Gr}\nolimits} 
\def\Fl{\mathop{\rm Fl}\nolimits} 
\def\Sym{\mathop{\rm Sym}\nolimits} 
\def\CAlg{\mathop{\rm CAlg}\nolimits} 
\def\End{\mathop{\rm End}\nolimits} 
\def\Map{\Maps} 

\def\Z{{\mathbf Z}} 
\def\Fp{{{\mathbf F}_p}} %
\def\Fq{{{\mathbf F}_q}} %
\def\Fqq{\overline {\mathbf F}_q} %
\def\Qlq{{\overline \Q_\ell}} %
\def\Q{{\mathbf Q}} 
\def\C{{\mathbf {C}}} 
\def\A{{\bf A}} 
\renewcommand{\P}[1][1]{\mathbf P^{#1}} 
\def\Gm{\mathbf {G}_\mathrm m} 
\newcommand{\twi}[1]{\langle #1\rangle} 

\def\Cc{{\mathcal{C}}}

\def\Ee{{\mathcal{E}}}
\def\Tt{{\mathcal{T}}}

\def\H{{\rm H}} 
\def\SH{\category{SH}} %
\def\HS{\category{HS}} 
\def\MHS{\category{MHS}} 
\def\pol{\category{pol}} 
\def\R{\mathrm{R}} 
\def\red{\mathrm{r}} 
\def\redx{(\red)} 
\def\DM{\category{DM}} 
\def\DMr{\DM_\red} 
\def\DK{\category{DK}} 
\def\DKr{\DK_\red} 
\def\DMrx{\DM_{\redx}} 
\def\DTM{\category{DTM}} 
\def\DTMr{\DTM_\red} 
\def\DTMrx{\DTM_{\redx}} 
\def\Par{\category{Par}} 
\def\DAR{\category{D}^{\mathrm{AR}}_{\mixed}} 
\def\Ld{\mathrm{L}} 

\def\MTM{\category{MTM}} 
\def\MTMr{\MTM_\red} 
\def\MTMrx{\MTM_{\redx}} 
\def\Tilt{\category{Tilt}} 
\def\ii{$\infty$}

\def\Gal{{\rm Gal}} 

\def\BR{\R_B} 

\def\Spec{\mathop{\rm Spec}} 
\newcommand{\Comp}{\category{Ch}} 
\newcommand{\M}{\mathrm{M}} 
\newcommand{\comp}{\mathrm{c}} 
\newcommand{\Frob}{\mathrm{Frob}} 
\def\bd{\mathrm{b}} 
\newcommand{\mix}{\mathrm{m}} 
\newcommand{\mixed}{\mathrm{mix}} 

\newcommand{\Tate}{\mathrm{T}} 
\def\E{{\rm E}} 
\newcommand{\an}{\mathrm{an}} 
\newcommand{\D}{\category{D}} 
\newcommand{\Ch}{\category{Ch}} 

\def\sbuildrel#1\over#2{\mathrel{\smash{\mathop{\kern0pt #2}\limits^{#1}}}}

\let\x\times
\let\ol\overline
\renewcommand{\t}{\otimes}

\renewcommand{\r}{\rightarrow}
\newcommand{\lr}{\longrightarrow}

\def\matrix#1{\null\,\vcenter{\normalbaselines
    \ialign{\hfil$##$\hfil&&\quad\hfil$##$\hfil\crcr
      \mathstrut\crcr\noalign{\kern-\baselineskip}
      #1\crcr\mathstrut\crcr\noalign{\kern-\baselineskip}}}\,}

\newdimen\harrowsize
\def\mapright#1{\smash{\mathop{\hbox to\harrowsize{\rightarrowfill}}\limits^{#1}}}

{\catcode`@=11
\gdef\cal{\fam\tw@}
\global\let\over\@@over
\global\let\atop\@@atop
\global\let\above\@@above
\global\let\overwithdelims\@@overwithdelims
\global\let\atopwithdelims\@@atopwithdelims
\global\let\abovewithdelims\@@abovewithdelims
\gdef\eqalign#1{\null\,\vcenter{\openup\jot\m@th
\ialign{\strut\hfil$\displaystyle{##}$&$\displaystyle{{}##}$\hfil
      \crcr#1\crcr}}\,}
\newskip\xcentering \global\xcentering=0pt plus 1000pt minus 1000pt
\gdef\eqalignno#1{\displ@y \tabskip\xcentering
  \halign to\displaywidth{\hfil$\@lign\displaystyle{##}$\tabskip\z@skip
    &$\@lign\displaystyle{{}##}$\hfil\tabskip\xcentering
    &\llap{$\@lign##$}\tabskip\z@skip\crcr
    #1\crcr}}
\global\def\cases#1{\left\{\,\vcenter{\normalbaselines\m@th
    \ialign{$##\hfil$&\quad##\hfil\crcr#1\crcr}}\right.}
\gdef\eqlabel#1{\refstepcounter{equation}\label{eqn--#1}\eqno\hbox{\@eqnnum}}
}

\def \nts#1{}

\def \journal#1
{ 
\noindent\colorbox{dunkelblau}{\parbox{\dimexpr\textwidth-2\fboxsep\relax}{#1}}
}

\begin{document}

\author{Jens Niklas Eberhardt and Jakob Scholbach}

\title{Integral Motivic Sheaves And Geometric Representation Theory}
\begin{abstract}
With representation-theoretic applications in mind, we construct a formalism of \emph{reduced motives} with integral coefficients. These are motivic sheaves from which the higher motivic cohomology of the base scheme has been removed. 
We show that reduced stratified Tate motives satisfy favorable properties including weight and t-structures.
We also prove that reduced motives on cellular (ind-)schemes unify various approaches to mixed sheaves in representation theory, such as Soergel--Wendt's semisimplified Hodge motives, Achar--Riche's complexes of parity sheaves, as well as Ho--Li's recent category of graded $\ell$-adic sheaves. 
\end{abstract}

\maketitle

\tableofcontents

\section{Introduction}
\subsection{Motivation}

Sheaves on manifolds are an important tool in geometric representation theory. For example, highest weight representations of a complex reductive Lie algebra can be described in terms of perverse sheaves on a flag manifold and representations of a reductive algebraic group correspond to equivariant perverse sheaves on an affine Grassmannian.

More refined formalisms of sheaves carry an additional notion of \emph{weights} and a \emph{Tate twist} functor that, in a very rough sense, provide an additional grading on the category. For example, mixed Hodge modules and mixed $\ell$-adic sheaves have a notion of weights via Hodge structures and eigenvalues of the Frobenius, respectively.
These have been put to great use in geometric representation theory. To name just two examples, the proof of the \emph{Kazhdan--Lusztig conjectures} crucially depends on the decomposition theorem for perverse sheaves for which weight considerations are essential while Beilinson--Ginzburg--Soergel's \emph{Koszul duality} for flag varieties relies on an additional grading provided by weights.

However, mixed Hodge modules and mixed $\ell$-adic sheaves have several drawbacks: first, they have characteristic zero coefficients and are hence not applicable in modular representation theory. Second, there are unwanted extensions between Tate objects which bear no representation-theoretic significance and yield technical problems.
Third, the choice of a fixed cohomology theory leaves open the question to which extent the resulting categories are independent of the coefficients or base field.

The goal of this paper is to introduce a formalism of \emph{mixed sheaves}, say $\D_{\mixed},$ that overcomes these problems. In particular, our proposed formalism works with (almost) arbitary coefficients, carries a six functor formalism and has no extension of Tate objects. Under appropriate assumptions, it is moreover independent of the base and specialises to the existing approaches to categories of mixed sheaves in the literature. 

\subsection{Reduced motives} A natural candidate for a category $\D_\mixed$ of mixed sheaves is the \emph{derived category of motivic sheaves} $\DM$ recalled in \refsect{motives}. It provides a universal home for cohomology and specialises to other formalisms of (mixed) sheaves. Building on that we define the category of \emph{reduced motives} in \refsect{DTMr.definition} as 
$$\DMr(X) := \DM(X) \t_{\DTM(S)} \grModL.$$
Here $X$ is a scheme or an ind-scheme of finite-type over a general base scheme $S$ such as $S=\Spec \Z$ and $\Lambda$ is any coefficient ring.
The above definition, which requires using \ii-categories, implements the idea that reduced motives are motives modulo the cohomology of the base $S.$ This results in independence of $S$ and removes unwanted extensions between Tate motives $\Lambda(n)$.
We prove in \refsect{DTMr.functoriality} that there is a six functor formalism for $\DMr$ compatible with the natural \emph{reduction functor} $\red : \DM \r \DMr.$

\subsection{Stratified Tate motives} For applications in geometric representation theory, we restrict our attention to more particular spaces and sheaves in \refsect{4}. We consider (ind-)schemes $X$ equipped with a well-behaved cellular stratification $\iota: X^+\to X$ into strata of the form $\A_S^{n} \x \GmX{S}^{m}$.  Our proposal for a formalism of mixed sheaves $\D_\mixed(X)$ on such spaces is the subcategory of \emph{reduced stratified Tate motives}
$$\DTMr(X,X^+)\subset \DMr(X)$$
which are roughly those motives that constant along the stratification. 
In particular, a reduced Tate motive on $\A^n_S$ is just a $\Z$-graded complex of $\Lambda$-modules (up to quasi-isomorphism):
$\DTMr(\A^n_S) = \grModL.$
Thus, reduced Tate motives have a strongly combinatorial flavor.

Reduced stratified Tate motives admit a perverse $t$-structure and a Chow weight structure. The heart of the $t$-structure, denoted by $\MTMr(X,X^+),$ can be considered an abelian category of \emph{mixed perverse sheaves}. The homotopy category of the heart of the weight structure, denoted by $\Ho(\DTMr(X, X^+)^{\wt=0}),$ is the additive category of \emph{pure complexes} and is, in many examples, generated by motives of resolutions of the stratum closures. We will show that both hearts completely determine $\DTMr(X, X^+).$ 

\theo (\thref{DTM.enoughtiltingperverse}, \thref{weight.complex.equivalence})
Under appropriate assumptions on the stratification,
there are equivalences of categories
\[\D^\bd(\MTMr(X,X^+)^{\comp}) \to \DTMr(X,X^+)^{\comp}\to \Ch^\bd(\Ho(\DTMr(X, X^+)^{\comp, \wt=0})).\]
\xtheo 

\subsection{Comparison Results}
In \refsect{comparison.results} we prove that reduced Tate motives refine and interpolate between various categories of mixed sheaves found in the literature. 

\theo \thlabel{DTM.comparisons}
Let $(X,X^+)$ be an (ind-)scheme over $S$ with a well-behaved cellular stratification. Depending on $S$ and the coefficient ring $\Lambda$, the category $\DTMr(X, X^+)$ is equivalent to the following categories:
\begin{enumerate}
\item
the category of (unreduced) stratified Tate motives $\DTM(X, X^+)$
in case $S = \Spec \mathbf{F}_{p^n}$ and $\Lambda = \Q$ or $\Lambda = \Fp$ (cf.~Soergel--Wendt \cite{SoergelWendt:Perverse} for $\Lambda=\Q$ and Eberhardt--Kelly \cite{EberhardtKelly:Mixed} for $\Lambda=\Fp$), see \thref{comparison.DM}.

\item
Achar--Riche's mixed derived category $\DAR(X^{\an})$ \cite{Achar:ModularPerverseSheavesII} which is the category of chain complexes of parity sheaves $\Par(X^{\an},X^+)$ \cite{Juteau:ParitySheaves2014} in case that $S = \Spec \C$, $\Lambda$ is a principal ideal domain and all strata of $X$ are affine spaces, see \thref{DTM.flat.paritysheaves}.

\item
Soergel--Wendt's category $\DTM_\calH(X,X^+)$ of semisimplified Hodge motives \cite{SoergelWendt:Perverse} in case $S = \Spec \C$ and $\Lambda=\C,$ see \thref{comparison.Hodge}.

\item 
Ho--Li's category of graded sheaves \cite{HoLi:Revisiting} (more precisely, the Tate objects therein) in case $S = \Spec \Fq$ and $\Lambda=\Qlq$, see \thref{comparison.graded}.
\end{enumerate}
\xtheo

Note that---unlike reduced (stratified Tate) motives---all the above theories are fixed to specific base schemes $S$ and, except for $\DAR$, also to a specific coefficient ring $\Lambda$. By contrast, under certain mild conditions, reduced stratified Tate motives over different base schemes $S$ are equivalent (\thref{independence.S}). This implies the equivalence of the above theories, whenever the coefficient ring agrees.

An advantage of reduced stratified motives over the category $\DAR$, which is the only one with integral coefficients so far, is the full-fledged six functor formalism. By comparison,  Achar--Riche's construction via chain complexes of parity sheaves offers only certain functors constructed by hand in \cite{Achar:ModularPerverseSheavesII}.

\subsection{Examples in Geometric Representation Theory}
Our formalism applies to various spaces used in geometric representation theory, such as (affine) partial flag varieties. 

Our running example is the flag variety $X=G/B$ of a split reductive group $G$ with Borel subgroup $B\subset G$ which has an affine stratifiction $X^+$ by $B$-orbits. 
Let $C=S/S^W$ be the coinvariant algebra, where $S=\Sym(X(T)_\Lambda)$ denote the symmetric algebra of the character lattice of the maximal torus $T\subset B.$ 
Assume that the torsion index of $G$ (for example $1$ in type $A_n,C_n$ or $2$ in type $B_n,D_n$) is invertible in $\Lambda$. Then the category of reduced stratified Tate motives is equivalent to the category of complexes of graded \emph{Soergel modules}
$$\DTMr(X,X^+)^\comp\cong\Ch^\bd(\SgrMod_C).$$
This immediatly follows from the comparison in \thref{DTM.flat.paritysheaves} and the corresponding proof of Soergel's \emph{Erweiterungssatz} in the literature, see \cite{Soergel:ICandRT,AcharRiche:ModularI}.
From this, we obtain the following diagram
\begin{center}
\begin{tikzcd}
  {\DTMr(X,X^+)_\C^\comp} \arrow[d,"\wr"] & {\DTMr(X,X^+)_\Z^\comp} \arrow[r] \arrow[l] & {\DTMr(X,X^+)_{\mathbf{F}_p}^\comp} \arrow[d,"\wr"] \\
 \D^\bd(\mathcal{O}^\Z(\mathfrak{g}^\Ld))                        &                                            &   \D^\bd(\mathcal{O}^\Z(G^\Ld))
  \end{tikzcd}
\end{center}
relating reduced stratified Tate motives to 
graded \emph{category $\mathcal{O}$} of the Langlands dual complex Lie algebra $\mathfrak{g}^\Ld/\C,$ see \cite{Soergel:Koinvarianten, Beilinson:Koszulduality} and graded \emph{modular category $\mathcal{O}$} of the Langlands dual group $G^\Ld/\mathbf{F}_p,$ see \cite{Soergel:ICandRT}. 

A very similar picture arises for reduced stratified Tate motives on the affine Grassmannian $\Gr_G$ with its stratification by Iwahori-orbits.
For $\Lambda=\C,$ it yields the derived graded principal block of finite dimensional representations of the quantum group $\operatorname{U}_q(\mathfrak{g}^\Ld/\C)$ at an odd root of unity, see \cite{Arkhipov:QuantumGroupsLoop2004}.
For $\Lambda=\Fp$, one obtains the derived graded principal block of algebraic representations of the algebraic group $G^\Ld/\mathbf{F}_p$ which amounts to a graded version of the \emph{Finkelberg–Mirkovi\'c conjecture}, see \cite{Achar:ParitySheavesAffine2015,Mautner:ExoticTiltingSheaves2018,Achar:ReductiveGroupsLoop2018}.

\subsection{Further directions}
\subsubsection{Realisation functor}
One conspicuous omission in the list of properties above is a realization functor, say for $X / \C,$
$$\DTMr(X, X^+) \stackrel {?} \r \D(X^\an, \Lambda).$$
Such a functor does exist for $\Lambda = \C$ (\thref{reduced.Hodge.realization}), but its existence for $\Lambda = \Z$ or $\Fp$ depends strongly on the way the strata in $X$ are glued together.
In upcoming work, we plan to investigate this topic as well as the related question how to ``unreduce'' reduced motives.
A concrete question in this direction is the following:

\ques
Let $X$ be the flag variety with the $B$-orbit stratification. Is there an equivalence
$$\DTM(X,X^+) \stackrel ? = \DTMr(X,X^+) \t_{\grModL} \DTM(S)?$$
\xques 
An affirmative answer to this question and similarly for affine flag varieties would seem to pave a way towards a solution of the Finkelberg–Mirkovi\'c conjecture via its graded version from \cite{Achar:ParitySheavesAffine2015,Mautner:ExoticTiltingSheaves2018,Achar:ReductiveGroupsLoop2018}, see \cite[Remark 2.13(3)]{Williamson:AlgebraicRepresentationsConstructible2017}.
\subsubsection{Equivariant ($K$-)motives} In a next step, we aim to extend our formalism to equivariant motives and $K$-motives. 
This would be a necessary first step for an integral motivic Satake correspondence, refining \cite{RicharzScholbach:Motivic}. It would also be useful in motivic Springer theory \cite{Eberhardt:MotivicSpringerTheory2022}. Moreover, it paves the way to study equivariant $K$-motives on flag varieties, the affine Grassmannian and the nilpotent cone which, conjecturally, yield an ungraded equivariant Koszul duality, see \cite{Eberhardt:KMotivesKoszulDuality2019}, a derived quantum $K$-theoretic Satake, see \cite{Cautis:QuantumKtheoreticGeometric2015}, and a $K$-theoretic motivic Springer theory.

For example, denote by $\mathcal{N}$ the nilpotent cone of a split reductive group $G.$ Then, certain reduced equivariant $K$-motives (motives) on $\mathcal{N}$ should yield the (graded) perfect derived category of the (graded) affine Hecke algebra $\mathbf{H}$ (resp. $\overline{\mathbf{H}}$).
\conj There are equivalences of categories
$$\DMr^{\mathrm{Spr}}(\mathcal{N}/(G\times \Gm))\cong \grPerf_{\overline{\mathbf{H}}} \text{ and }\DKr^{\mathrm{Spr}}(\mathcal{N}/(G\times \Gm))\cong \Perf_{\mathbf{H}}.$$
\xconj
In light of \thref{comparison.DM}, the conjecture is a refinement of \cite[Theorem 1.3]{Eberhardt:SpringerMotives2021} where a similar statement is shown for $\DM^{\mathrm{Spr}}(\mathcal{N}_\Fp/(G\times \Gm))_\Q.$

\subsection*{Acknowledgements.} We thank Vladimir Sosnilo, 
Matthias Wendt,
and Yifei Zhao 
for helpful discussions. 
J.N.E. was supported by Deutsche Forschungsgemeinschaft (DFG), project number 45744154,  Equivariant K-motives and Koszul duality. 
J.S. was supported by Deutsche Forschungsgemeinschaft (DFG), EXC 2044–390685587, Mathematik Münster: Dynamik–Geometrie–Struktur. 

\section{Recollections}

In this section we recall some notions related to \ii-categories, as well as a coherent formulation for the six functor formalism for motives, and a description of Tate motives as modules over a graded $\E_\infty$-ring spectrum.

\subsection{Categorical generalities}

\subsubsection{\ii-categories}
\label{sect--ii.categories}

Concerning \ii-categories, we use the standard terminology of \cite{Lurie:Higher, Lurie:HA}. Thus, $\PrSt$ denotes the \ii-category of stable presentable \ii-categories and colimit-preserving functors. This category is endowed with the Lurie tensor product \cite[4.8.2.10,~4.8.2.18]{Lurie:HA}, with the monoidal unit being the category $\Sp$ of spectra. 
For any $C \in \PrSt$, and $c, d \in C$, there is a mapping spectrum of maps from $c$ to $d$, denoted by $\Maps_C(c, d)$. 
We denote the homotopy category of an \ii-category $C$ by $\Ho(C)$.
The hom-sets in this category are denoted by $\Hom_{\Ho(C)}$ or just $\Hom_C$.

Throughout we write $\Mod_A(C)$ for the \ii-category of $A$-modules, for any commutative algebra object $A \in \CAlg(C)$, for any symmetric monoidal \ii-category $C$. 
For $A \in \CAlg(\PrSt)$, the category $\Mod_A(\PrSt)$ is a symmetric monoidal \ii-category with tensor product denoted by $- \t_A -$ \cite[4.5.2.1]{Lurie:HA}.
The \ii-category of chain complexes of $\Lambda$-modules (for a commutative ring $\Lambda$) up to quasi-isomorphism is denoted $\Mod_\Lambda$ \cite[1.3.5.8]{Lurie:HA}.
The homotopy category $\Ho(\ModL)$ is the classical unbounded derived category of $\Lambda$-modules.
For $C \in \PrSt$, we denote the subcategory consisting of compact objects by $C^\comp$.
For example, $\Mod_\Lambda^\comp$ is the \ii-category of perfect complexes of $\Lambda$-modules, denoted $\Perf_\Lambda$.
By default, all functors are derived. This applies in particular to the tensor product, so that expressions such as $M \t_\Lambda N$ denote the derived tensor product, even if $M$ and $N$ are ordinary $\Lambda$-modules.

We write $\PrStL := \Mod_{\Mod_\Lambda}(\PrSt)$ for the \ii-category of presentable stable $\Lambda$-linear categories.

Given a commutative monoid $A \in \CAlg(\PrSt)$, any (stable presentable) \ii-category acted upon by $A$, i.e., any $C \in \Mod_A(\PrSt)$ is canonically enriched over $A$: we define the enriched mapping object $\IMaps_C(c, d) \in A$ to be the object (in $\Mod_A$) representing the functor $A^\opp \r \Sp$, $a \mapsto \Maps_C(a \t c, d)$, i.e., satisfying
$$\Maps_A(a, \IMaps_C(c, d)) = \Maps_C(a \t c, d).$$
\nts{Here we use the \cite[Proposition~5.5.2.2]{Lurie:Higher} to ensure that this functor is representable: $a \mapsto a \t c$ preserves colimits, hence the fuctor $A^\opp \r \Sp$ preserves limits.} 

\lemm
\thlabel{adjunctions.tensor}
Let $A \in \CAlg(\PrSt)$ be \emph{rigid}. Let
$$L : C \rightleftarrows C' : R$$
be an adjunction with $L$ being a map in $\Mod_A(\PrSt)$, i.e., an $A$-linear colimit-preserving functor. Assume that $R$ also preserves colimits.
\nts{Equivalently, if $C$ and $C'$ are compactly generated, $L$ preserves compact objects.}
Then $R$, which a priori is only a lax $A$-linear functor, is in fact $A$-linear.
Moreover, for any $D \in \Mod_A(\PrSt)$, there is an adjunction
$$L \t_A \id_D : C \t_A D \rightleftarrows C' \t_A D : R \t_A \id_D.$$
\xlemm

\pf
The first claim is \cite[Chapter I, Lemma~9.3.6]{GaitsgoryRozenblyum:StudyI} or \cite[Lemma~3.5]{BenZviNadler:Character}. 
Thus the expression $R \t_A \id_D$ makes sense to begin with.
The claim about the adjunction $L \t \id_D \dashv R \t \id_D$ holds by the characterization of adjunctions in terms of the triangle identities \cite[Digression~2.1.2]{RiehlVerity:Elements}.
\nts{In the absence of $A$ being rigid, we can also alternatively assume that the underlying functors in $\Cat_\infty$ are adjoint and such that the unit and counit maps $\id \r RL$, $LR \r \id$ are $A$-linear.}
\xpf

Recall that the \emph{lax limit} of a functor $f : C \r C'$ in $\PrSt$, denoted $\laxlim (C \stackrel f \r C')$, can be defined as the pullback of the following diagram:
$$\xymatrix{
& \Fun(\Delta^1, C') \ar[d]^{\ev_1} \\
C \ar[r]^{f} & C'.
}\eqlabel{laxlim}$$
Thus, an object in $\laxlim f$ is a triple $(c \in C, c_1' \r c_2' \in C', \alpha: f(c) \cong c_2')$. Under an equivalence of categories, objects in $\laxlim f$ can be described as triples $(c, c', c' \r f(c))$, where the map is arbitrary (not necessarily an isomorphism).

\lemm
\thlabel{laxlim.tensor}
Let $f : C \r C'$ be a map in $\Mod_A(\PrSt)$, for some commutative algebra object $A \in \CAlg(\PrSt)$.
Then $\laxlim f$ is naturally also an object in $\Mod_A(\PrSt)$.
Moreover, if $A$ is rigid and $B$ is an $A$-module that is compactly generated (more generally, dualizable in $\PrSt$, i.e., disregarding the $A$-action), then 
$$(\laxlim f) \t_A B = \laxlim (C \t_A B \stackrel {f \t \id_B} \r C' \t_A B).$$
\xlemm

\pf
The $A$-module structure on $\laxlim f$ arises since both maps in \refeq{laxlim} are $A$-linear, and since the forgetful functor $\Mod_A(\PrSt) \r \PrSt$ preserves limits.
We now use the following generalities \cite[Chapter~1, 4.1.6, 7.3.2, 9.4.4]{GaitsgoryRozenblyum:StudyI}: 
if $A$ is rigid, dualizability in $\Mod_A(\PrSt)$ is equivalent to being dualizable in $\PrSt$.
In addition, any compactly generated category is dualizable in $\PrSt$.
By dualizability of $B$ (over $A$) we also have $\Fun(\Delta^1, C') \t_A B = \Fun(\Delta^1, C' \t_A B)$.
Then, use that tensoring with dualizable objects preserves limits since $- \t_A B = \Fun_A(B^\dual, -)$, i.e., tensoring with $B$ is equivalent to considering the category of $A$-linear functors (within $\PrSt$) out 
of the ($A$-linear) dual of $B$.
\xpf

\subsubsection{Graded objects}

\defi
For a stable \ii-category $C$, write $\gr C$ for the category of $\Z$-graded objects in $C$, i.e, $\gr C := \Fun(\Z, C)$, where here $\Z$ is regarded as a discrete category.
We write
$$(-)_r := \ev_r : \gr C \r C$$
for evaluation at graded degree $r$, i.e., precomposition with $\{r \} \r \Z$.
We also let 
$$\twi r : \gr C \r \gr C$$
be the precomposition with $\Z \r \Z$, $m \mapsto m+r$.
Thus it shifts the grading by $r$, i.e., $(X \twi r)_m = X_{m+r}$.
\xdefi

\rema
The functor $\ev_r$ has a right and a left adjoint. These two adjoints agree (since $C$ is pointed) and are given by $C \r \gr C, X \mapsto (\dots, 0, X, 0, \dots)$, (insert $X$ in degree $r$ and the zero object elsewhere).

We will usually not distinguish between an object $X \in C$ and the same object, regarded as being concentrated in graded degree 0.
For example, for $X \in C$, the above graded object will be denoted by $X \twi {-r}$.
\xrema

\rema
If $C$ is symmetric monoidal, then so is $\gr C$ by means of the Day convolution product, with respect to the symmetric monoidal structure on $\Z$ given by addition.
In particular $\grModL := \gr(\ModL)$ is a commutative algebra object in $\PrSt$.
We will denote by $\grPrStL := \Mod_{\grModL}(\PrSt)$ its category of modules, i.e., the \ii-category of presentable stable $\Lambda$-linear and $\Z$-graded \ii-categories (with functors preserving these structures).

The functor $\twi 0 : C \r \gr C$ is symmetric monoidal.
In particular, the monoidal unit of $\gr C$ is given by $1\twi 0$.
Thus, its (right and left) adjoint $\ev_0$ is symmetric lax monoidal (and symmetric oplax monoidal, but not symmetric monoidal: $\ev_0((A_n) \t (B_n)) = \bigoplus_{a+b=0} A_a \t B_b \ne A_0 \t B_0$.)
\xrema

\subsubsection{t- and weight structures}

We briefly recall \emph{weight structures} and \emph{t-structures}, since the category of reduced stratified Tate motives enjoy both these structures.
The definitions are very similar, so let $\epsilon := +1$ in case of a t-structure, and $\epsilon := -1$ for a weight structure.
A \emph{t-structure} (resp., a \emph{weight structure}) on a stable \ii-category $C$ is a pair of full idempotent closed subcategories $(C^{\ge 0}, C^{\le 0})$ such that the following conditions hold. 
\begin{enumerate}
\item 
\label{item--shifts}
for each object $c \in C$ there is a fiber sequence, with $c^{\le 0} \in C^{\le 0}$, $c^{\ge 0} \in C^{\ge 1}$,
$$c^{\le 0} \r c \r c^{\ge 0}[-\epsilon],$$

\item
The subcategories $C^{\le 0}$ (resp.~$C^{ \ge 0}$) are stable under shifts by$[\epsilon]$ (resp.~$[-\epsilon]$).
(Thus, we are using a cohomological convention for t-structures and a homological for weight structures.)

\item
for $c^{\le 0} \in C^{ \le 0}$ and $c^{\ge 0} \in C^{\ge 0}$ the mapping spectra satisfy
$$\Hom_{\Ho(C)}(c^{\le 0}, c^{\ge0}[-\epsilon]) = \pi_0\Map_C(c^{\le 0}, c^{\ge 0}[-\epsilon]) = 0.$$
\end{enumerate}

In the presence of \refit{shifts}, the last condition is equivalent to $$\Hom_{\Ho(C)}(c^{\le 0}, c^{\ge 0}[-n]) = 0\text{ for all }n \in \Z\text{ with }\epsilon n > 0.$$

Weight structures on \ii-categories were introduced in \cite{Sosnilo:Theorem}; the notion is due to Bondarko \cite{Bondarko:Weight} and Pauksztello \cite{Pauksztello:CoTStructures} in the context of triangulated categories. Note that we use the cohomological convention for t-structures and homological convention for weight structures here.

The \emph{heart} of a weight-structure or t-structure is the full subcategory $C^{= 0} := C^{\le 0} \cap C^{\ge 0} \subset C$.
It is an additive \ii-category for weight structures, and an abelian category for t-structures.
We indicate the aisles and the heart of a weight structure also by $C^{\wt \le 0}$, $C^{\wt = 0}$, and the ones of a t-structure by $C^{\tt \le 0}$, $C^{\tt = 0}$ etc. 

A functor between two such categories with weight (or t-)structures is \emph{weight- (or t-)exact} if it preserves the two given subcategories.

For a bounded weight structure $\wt$ on a stable \ii-category $C$ (i.e., such that $\bigcup_{n \in \Z} C^{\wt \ge 0}[n] = C = \bigcup_{n \in \Z} C^{\le 0}[n]$), there is an essentially unique exact functor, called \emph{weight complex functor},
$$C \r \Comp^\bd(\Ho(C^{\wt = 0}))$$
that restricts to the identity on $C^{\wt = 0}$ \cite[Corollary~3.0.4]{Sosnilo:Theorem}. Here the target category is the \ii-category underlying the dg-category formed by bounded complexes in the additive category $\Ho(C^{\wt = 0})$.

Similarly, if $(C^{\tt \le 0}, C^{\tt \ge 0})$ is a t-structure, there is a unique exact functor (up to equivalence), called \emph{realization functor},
$$\D^\bd(C^{\tt = 0}) \r C\eqlabel{realization.functor}$$
whose restriction to $C^{\tt = 0}$ is the identity \cite[Remark~7.60]{BunkeCisinskiKasprowskiWinges:Controlled}.

Weight structures can be extended to ind-completed categories by means of the following lemma, which is due to Bondarko (in the setting of triangulated categories; the translation to stable \ii-categories is routine).

\lemm \cite[Theorem~4.1.2]{Bondarko:Morphisms}
\thlabel{weights.Ind}
Let $C$ be an essentially small stable \ii-category with a weight structure.
\begin{enumerate}
\item 
\label{item--Ind.C.weights.description}
Its Ind-completion $\Ind C$ (cf.~\cite[Proposition~1.1.3.6]{Lurie:HA}) carries a weight structure such that $(\Ind C)^{\wt \ge 0}$ (resp.~$(\Ind C)^{\wt \le 0}$) is the smallest full subcategory containing $C^{\wt = 0}$, and stable under coproducts, extensions and shifts $[+1]$ (resp. $[-1]$).
\item
The heart $(\Ind C)^{\wt = 0}$ is the full subcategory of $\Ind(C^{\wt = 0})$ containing $C^{\wt = 0}$ and arbitrary coproducts.
\item
A functor $\Ind C \r D$ taking values in a stable \ii-category with weight structure is weight exact iff its restriction to $C$ is so.
\end{enumerate} 
\xlemm

\exam
\thlabel{ModL.t.weight.structure}
The category $\ModL = \Ind(\Perf_\Lambda)$ has its natural t-structure whose heart is the usual category of $\Lambda$-modules. The t-structure restricts to one on $\Perf_\Lambda$ iff $\Lambda$ is a regular coherent ring \cite[Proposition~6.6]{HemoRicharzScholbach:Constructible}, for example a regular Noetherian ring. 
The category $\Perf_\Lambda$ also has a weight structure whose heart is the (additive) category of finitely generated projective $\Lambda$-modules \cite[Example~3.1.6]{Sosnilo:Theorem}. It gives a weight structure on $\ModL$ by ind-extension as in \thref{weights.Ind}.

The category $\grModL$ has a t-structure and weight-structure such that the evaluation functors $\ev_n : \grModL \r \ModL$ are t- and weight exact.
(Thus, the graded degree does not affect the t- or weight degree of an object.)   
\xexam

\subsection{Motivic sheaves}

\label{sect--motives}

\conv
\thlabel{convention.S}
Throughout the entire paper, we fix a scheme $S$ that is supposed to be connected, smooth and of finite type over  $\Spec O$, where $O$ is a field or a Dedekind ring.
(Thus, $S$ is itself regular, Noetherian and of finite Krull dimension; in practice we only care about the case $S = \Spec O$ itself.)
The category $\Sch_S^\ft$ is the category of $S$-schemes of finite type. We refer to objects in this category just as $S$-schemes or even just as schemes.

We also fix a regular coherent coefficient ring $\Lambda$.
\xconv

\defi
For any scheme $X$, the category of \emph{motives over $X$} is defined as
$$\DM(X) := \DM(X)_\Lambda := \Mod_{\M\Lambda}(\SH(X)),$$
the category of modules, in the stable $\A^1$-homotopy category, over the motivic ring spectrum $\M\Lambda$ representing motivic cohomology (with coefficients in a commutative ring $\Lambda$).
\xdefi

Building on the seminal work of Morel--Voevodsky, the category $\SH$ was developed by Ayoub \cite{Ayoub:Six1}.
For schemes over a Dedekind ring, the ring spectrum $\M\Z$ was introduced by Spitzweck \cite{Spitzweck:Commutative}.
For any commutative ring $\Lambda$ (we will mostly use $\Z$ and fields), this gives a ring spectrum $\M\Lambda$ by scalar extension. 
A key asset of these categories is the six functors formalism, which has been developed by Ayoub, Cisinski--Déglise \cite{CisinskiDeglise:Triangulated},
and in an \ii-categorical form by Khan \cite{Khan:Motivic}, building upon work of Gaitsgory--Rozenblyum \cite{GaitsgoryRozenblyum:StudyI}.
In \cite[Appendix~A]{RicharzScholbach:Motivic}, this was extended to allow for a coherent handling of monoidal aspects, including projection formulas.
We briefly recall this extension (which just adds an $\epsilon$ to the previously existing literature):
let $\Sch_S^{\ft, \x} \r \Fin_*$ be the symmetric monoidal category associated to with the cartesian monoidal structure on $\Sch_S^\ft$.
Let $\Sch_S^{\ft, \x, \dual} \r \Fin_*^\opp$ be the associated dual fibration. The opposite of that map encodes the usual symmetric monoidal structure on $(\Sch_S^{\ft})^\opp$. 
In order to encode *- and !-functoriality at the same time, one uses the category of correspondences. We refer to \cite{GaitsgoryRozenblyum:StudyI} for a full discussion of this category.
Here, we only point out that $\Corr := \Corr(\Sch_S^{\ft, \times, \dual})^\proper_{\sep, \all}$ is an $(\infty,2)$-category whose objects are the objects in $\Sch_S^{\ft, \times}$ (i.e., sequences of objects $X := (X_1, \dots, X_n)$ with $X_i \in \Sch_S^\ft$).
In this category, 1-morphisms from $X$ to $Y$ are spans $Y \stackrel g \gets Z \stackrel f \r X$, with $f = (f_i : Z_i \r X_i)$ being a collection of separated maps (i.e., the image of $f$ in $\Fin_*$ is an identity map) and $g$ being an arbitrary map. 
In $\Corr$, 2-morphisms between two 1-morphisms $Y \gets Z \r X$ and $Y \gets Z' \r X$ are maps $Z \stackrel z \r Z'$ (fitting into the obvious commutative diagrams) that map to an identity in $\Fin_*$, and whose components $Z_i \stackrel{z_i} \r Z'_i$ are proper.
There is a symmetric monoidal structure on $\Corr$ which on the level of objects is given by concatenating sequences of objects.

Then there is a lax symmetric monoidal functor 
$$\DM^*_! : \Corr := \Corr(\Sch_S^{\ft, \times, \dual})^\proper_{\sep, \all} \r \PrStL$$
whose value at some scheme $X$ is the category $\DM(X)$ mentioned before.

\rema
\thlabel{Corr.explained}
Among other things, this functor encodes the following data.
\begin{itemize}
\item 
For a separated map $f$ and any map $g$, there are two composable morphisms in $\Corr$:
$$\xymatrix{
X \x_Y Z \ar[r]^{g'} \ar[d]^{f'} & X \ar[d]^f \ar@{=}[r] & X \\
Z \ar[r]^g \ar@{=}[d] & Y 
\\
Z
}$$
The evaluation of $\DM^*_!$ at the top right correspondence is a functor $f^* : \DM(Y) \r \DM(X)$, the one for the lower correspondence is a functor $g_! : \DM(Z) \r \DM(Y)$.
The functoriality of $\DM^*_!$ thus encodes the base-change formula
$$f^* g_! = g'_! f'^*.$$
  \item

  Any $S$-scheme $Y$ is a commutative comonoid in $\Sch_S^{\ft}$ (since the cartesian structure is used), with comultiplication given by the diagonal $Y \stackrel{\Delta} \r Y \x_S Y$. 
  The lax monoidality of $\DM^*_!$ then yields a functor
$$\DM(Y) \t_{\ModL} \DM(Y) \stackrel \boxtimes \r \DM(Y \x_S Y).$$
Appending $\Delta^*$, we obtain that $\DM(Y)$ becomes a symmetric monoidal \ii-category.
The monoidal unit in $\DM(Y)$ is denoted $\Lambda$ (recall that $\Lambda$ is the coefficient ring).

  \item
  Any map $f : X \r Y$ in $\Sch_S^\ft$ turns $X$ into a $Y$-comodule,
  with coaction given by $X \stackrel \Delta \r X \x_S X \stackrel{f \x \id} \r X \x_S Y$.
  The lax monoidality of $\DM^*_!$ implies that $\DM(Y)$ is a $\DM(X)$-module by means of $f^*$.
With these structures, $f$ is a map of $Y$-coalgebras so that the evaluation of the lax symmetric monoidal functor $\DM^*_!$ encodes the projection formula
  $$f_! (A \t f^* B) = f_! A \t B \ \ (A \in \DM(X), B \in \DM(Y)).$$
\end{itemize}
\xrema

\rema
\thlabel{DM.further.properties}
In addition to the functoriality encoded by correspondences, $\DM$ satisfies the following properties:
\begin{enumerate}
\item \emph{Homotopy invariance}: for the projection $p : \A^1_X \r X$, the following functor is fully faithful:
$$p^* : \DM(X) \r \DM(\A^1_X).$$

\item
\emph{Tate twists}: for the projection $q : \GmX X \r X$, the object $\Lambda(1) := \fib (q_! q^! \Lambda \r \Lambda)[1]$ is $\t$-invertible with dual denoted $\Lambda(-1)$. We put $\Lambda(n) := \Lambda(1)^{\t n}$ for $n \in \Z$.

\item
\label{item--localization}
\emph{Localization}: for a closed immersion $i: Z \r X$ with complement $j: U \r X$, the (co)units of the adjunctions above assemble into so-called \emph{localization} homotopy fiber sequences
$$i_! i^! \r \id \r j_* j^*,\eqlabel{localizationOne}$$
$$j_! j^! \r \id \r i_* i^*.\eqlabel{localization}$$

\item
\label{item--algebraic.cycles}
\emph{Algebraic cycles}: for $X$ smooth over $S$, the Hom-groups are given by higher Chow groups of Bloch (extended to schemes over Dedekind rings by Levine):
$$\Hom_{\Ho(\DM(X))}(\Lambda, \Lambda(n)[m]) =: \H^m(X, \Lambda(n)) = \CH^n(X, 2n-m)_\Lambda.$$
For later use, we note that this group vanishes for $m>2n$.
\end{enumerate}
\xrema

For the purpose of defining reduced motives, it will be useful to uniformly keep track of the presence of the $\DM(S)$-action on all categories of motives.
Indeed, the category $\Sch_S^{\ft}$ identifies with $\Mod_S(\Sch_S^\ft)$, since $S$ is the monoidal unit. Applying the lax symmetric monoidal functor $\DM$, we conclude the existence of a functor
$$\DM^*_! : \Corr \r \Mod_{\DM(S)}(\PrStL).$$

\subsection{Motives on ind-schemes}
\label{sect--DM.ind-schemes}

The functor $\DM_! := \DM^*_!|_{\Sch_S^{\ft, \sep}}$ (schemes with only separated maps) can be used to define motives on \emph{ind-schemes} by a left Kan extension \cite[\S2.3]{RicharzScholbach:Intersection}: 
$$\xymatrix{
(\Sch_S^{\ft})^\sep \ar[d] \ar[r]^{\DM_!} & \Mod_{\DM(S)}(\PrSt) \\
(\IndSch_S^{\ft})^\sep \ar@{.>}[ur]^{\DM_!}
}\eqlabel{DM.indschemes}$$
In other words, for an ind-scheme $X = \colim X_i$, the category of motives is given by 
$$\DM(X) = \colim \DM(X_i),$$
where the transition functors in this colimit are !-pushforwards along the closed embeddings $X_i \r X_j$.
By definition, $\DM(X)$ is a presentable stable \ii-category, and a module over $\DM(S)$.
The subcategory $\DM(X)^\comp$ of compact objects can be thought of as the union of the categories $\DM(X_i)^\comp$, again using (the fully faithful) !-pushforwards to form the union.

For a map $f : X \r Y$ between ind-schemes, there is always an adjunction $f_! : \DM(Y) \rightleftarrows \DM(X) : f^!$, while the (adjoint) functors $f^*$ and $f_*$ only exist for schematic maps, see \cite[Theorem~2.4.2]{RicharzScholbach:Intersection} for further details.

\subsection{Tate motives as modules over a graded $\E_\infty$-algebra}

In this section, we recall the description of the category of Tate motives as a category of modules.
This presentation will be crucial for the definition of reduced motives.
Everything in this section is due to Spitzweck \cite{Spitzweck:Mixed}; we include certain proofs for the convenience of the reader.

\defi
For a scheme $X$, the category $\DTM(X)$ of \emph{Tate motives} is defined to be the presentable subcategory of $\DM(X)$ generated by the $\Lambda(n)$ for $n \in \Z$.
\xdefi

\lemm (\cite[Proposition 4.2]{Spitzweck:Mixed}, \cite[\S8]{Spitzweck:Commutative})
There is a commutative monoid in $\gr\DM(S)$, denoted $\PM\Lambda$, whose underlying object in $\gr\DM(S)$ is the one whose component in graded degree $r$ is $\Lambda(r)$.
\xlemm
Thus, morally speaking $\PM\Lambda = \bigoplus_{n \in \Z} \Lambda(n)$. 
Using the notation of \refsect{ii.categories}, there is an adjunction
$$L: \grModL \rightleftarrows \gr\DM(S) : R,$$
whose left adjoint $L$ satisfies $\Lambda \mapsto \Lambda$. 
Here we use that $S$ is connected so that $\End_{\DM(S)}(\Lambda, \Lambda) = \Lambda$.
Since $R$ is lax symmetric monoidal, it preserves commutative monoids, so that
$$A := R(\PM\Lambda)\eqlabel{A}$$
is a commutative monoid in the \ii-category $\grModL$.
Roughly, one can think of it as a $\Z$-graded commutative differential graded $\Lambda$-algebra. (Unless $\Q \subset \Lambda$, this is only a rough analogy since it is not usually possible to strictify a commutative algebra in the \ii-category $\Mod_\Lambda$ to a commutative dg-algebra.)

We compute $A$ as follows: by adjunctions we have
$$\eqalign{
A_r & = \Map_{\ModL}(\Lambda, A_r) \cr
& = \Map_{\grModL}(\Lambda \twi {-r}, A) \cr
& = \Map_{\gr\DM(S)}(\Lambda \twi {-r}, \PM\Lambda) \cr
& = \Map_{\DM(S)}(\Lambda, (\PM\Lambda)_r) \cr
& = \Map_{\DM(S)}(\Lambda, \Lambda(r)).}\eqlabel{compute.A}$$
Thus the $r$-th graded component $A_r$ is a chain complex whose $n$-th cohomology is $\H^n(S, \Lambda(r))$.

By the universal property of the category of modules \cite[\S3.3.3]{Lurie:HA}\nts{See \cite[\S3.3.9]{Robalo:Theorie} for pointers and further explanation}, the functor $R$ induces a limit-preserving, accessible functor, again denoted $R$, which by the adjoint functor theorem admits a left adjoint:
$$\tilde L : \Mod_A(\grModL) \rightleftarrows \Mod_{\PM\Lambda}(\gr\DM(S)) : R.$$
The left adjoint $\tilde L$ satisfies $\tilde L (A \t V) = \PM\Lambda \t L(V)$ for $V \in \grModL$.
In particular, $\tilde L$ is symmetric monoidal, so that $\tilde L (A) = \PM\Lambda$.
Henceforth we abbreviate 
$$\Mod_A := \Mod_A(\grModL).$$
(If $A$ can be represented by a $\Z$-graded commutative differential graded $\Lambda$-algebra, the homotopy category of $\Mod_A$ is the category of graded complexes with an $A$-module structure, up to quasi-isomorphism.)

\lemm (\cite[Theorem~4.5]{Spitzweck:Mixed}, \cite[Corollary~8.3]{Spitzweck:Commutative})
The composite 
$$F := \ev_0 \circ \tilde L : \Mod_A \stackrel {\tilde L} \r \Mod_{\PM\Lambda}(\gr\DM(S)) \stackrel{\ev_0} \r \DM(S)$$
induces an equivalence of symmetric monoidal \ii-categories
$$\Mod_A \cong \DTM(S), A \twi r \mapsto \Lambda(r).\eqlabel{Mod.A.DTM}$$
\xlemm

\pf
The category $\grModL$ is compactly generated by the objects $\Lambda \twi r$, $r \in \Z$.
\nts{Generally $\gr C$ is compactly generated by compact objects in $C$, put in arbitrary graded degrees.}
The category $\Mod_A(\grModL)$ is then compactly generated by $A \twi r = A \t \Lambda \twi r$.
Our functor sends $A \twi r$ to $\ev_0(\PM\Lambda \twi r) = \ev_{r}(\PM\Lambda) = \Lambda(r)$, which is a compact object in $\DM(S)$. 

Any functor in $\PrSt_\omega$ (compactly generated presentable stable \ii-categories and continuous functors preserving compact objects) is fully faithful iff its restriction to compact objects is fully faithful.
Thus, in our case it is enough to check full faithfulness of $F$ restricted to our family of generators, $A \twi r$.
That is, we have to ensure that the mapping spectra
$$\Map_{\Mod_A(\grModL)}(A \twi s, A \twi r) \r \Map_{\DM(S)}(\Lambda(s), \Lambda(r))$$
are isomorphic. Indeed, the left hand side is just $\Map_{\ModL}(\Lambda, (A \twi {r-s})_0) = A_{s-r}$. This is precisely the right hand mapping space by design of $A$, cf.~\refeq{compute.A}.

Given the full faithfulness of $F$, the we obtain an equivalence of \ii-categories since the generators of $\DTM(S)$, $\Lambda(n) = F(A \twi{n})$ are in the image.
Being a composite of the symmetric monoidal functor $\tilde L$ and the lax symmetric monoidal $\ev_0$, $F$ is lax symmetric monoidal.
It remains to observe that the lax structural maps $F(M \t_A M') \r F(M) \t_\Lambda F(M')$ are isomorphisms. Since both sides are colimit-preserving in $M$ and $M'$, it suffices to check this for generators of the form $M = A \twi r$, $M' = A \twi {r'}$, where it is clear.
\xpf

\defi
A commutative monoid object $A$ in $\grModL$ is called of \emph{Tate type} if
\begin{enumerate}

\item 
the unit map
$$\Lambda \r A$$
induces an isomorphism after applying $\ev_0$, i.e., $\Lambda \r A_0$ is an isomorphism (in $\ModL$, i.e., a quasi-isomorphism of chain complexes), and

\item
the evaluations $A_r = 0$ for $r < 0$.
\end{enumerate} 
\xdefi

\lemm
The algebra $A := R(\PM\Lambda)$ is of Tate type.
\xlemm

\pf
The $n$-th cohomology of the complex
$$A_r = \Map_{\DM(S)}(\Lambda, \Lambda(r))$$
is isomorphic \cite[Corollary~7.19]{Spitzweck:Commutative} to
$$\H^n(S_\Zar, \calM(r)),$$
where $\calM(r)$ is the complex that has in cohomological degree $i$ the cycles $z^r(-, 2r-i) \t_\Z \Lambda$.
For $r < 0$, this complex is defined to be zero, and for $r = 0$, this complex is isomorphic to $\Lambda$, whose cohomology is just $\Lambda^{\pi_0(S)} \stackrel{\text{\ref{convention.S}}} = \Lambda$.
\xpf

\lemm \cite[Lemma~4.10]{Spitzweck:Mixed}
\thlabel{augmentation.map}
For an algebra of Tate type $A$, the map
$$A \r A_0 \cong \Lambda  (:= \Lambda \twi 0)$$
is part of a natural map of commutative monoid objects in $\grModL$.
We call this map the \emph{augmentation map}.
\xlemm

\pf
The inclusion $C \r \gr C$ in graded degree 0 is symmetric monoidal, so that its adjoint $\ev_0$ is lax symmetric monoidal: we have $\ev_0(X \t Y) = \bigoplus_{i+j=0} \ev_i X \t \ev_j Y \r \ev_0X \t \ev_0 Y$.
The restriction of $\ev_0$ to the subcategory of graded objects $X$ satisfying $\ev_r X = 0$ for $r < 0$, 
is a symmetric monoidal functor since for such $X$ the above maps are isomorphisms.
\xpf

\section{Reduced motives}

In this section, we define reduced motives and establish their basic functoriality.
The general idea of reduced motives is to suppress all motivic cohomology coming from the base scheme $S$ (i.e., the one that is present in $A$ in \refeq{A}), but leave the remainder of the motivic formalism intact. 

We continue to fix a base scheme $S$ and a coefficient ring $\Lambda$ as in \thref{convention.S}.

\subsection{Definition and immediate properties}
\label{sect--DTMr.definition}

For any $S$-scheme $X$, the (stable \ii-)category $\DM(X)$ of motives on $X$ is a module over $\DM(S)$ by means of the pullback $f^*$ along the structural map $f : X \r S$.
By restriction, $\DM(X)$ and also its full subcategory $\DTM(X)$, become $\DTM(S)$-modules.
This module structure is compatible with *-pullback (resp.~!-pushforward) along arbitrary (resp.~separated) maps and therefore continues to exist for $X$ being an ind-scheme (cf.~\refsect{motives}). 

On the other hand, using the augmentation map $A \r \Lambda$ (\thref{augmentation.map}), we also have the following functor (which one can think of as modding out the augmentation ideal, but the functor is derived):
$$\DTM(S) \stackrel{\text{\refeq{Mod.A.DTM}}} \cong \Mod_A(\grModL) \stackrel{\Lambda \t_A -} \r \Mod_{\Lambda}(\grModL) = \grModL.$$

\defi
Let $X / S$ be a scheme or an ind-scheme.
The category of \emph{reduced motives} on $X$ is defined as
$$\DMr(X) := \DMr(X)_\Lambda := \DM(X)_\Lambda \t_{\DTM(S)_\Lambda} \grModL.\eqlabel{DMr.definition}$$
Here the tensor product is formed in $\Mod_{\DTM(S)}(\PrSt)$, cf.~\refsect{ii.categories} for notation.
In the same vein, the category of \emph{reduced Tate motives} is defined as
$$\DTMr(X) := \DTM(X) \t_{\DTM(S)} \grModL.$$
\xdefi

\rema \phantom{a}
\thlabel{DTM.flat.examples}
\begin{enumerate} [wide, labelwidth=!]
\item
In order to form the above tensor product, it is crucial to use an \ii-categorical enhancement, as opposed to a mere triangulated category structure on $\DM(X)$. 

\item 
By definition,
$$\DTMr(S) = \grModL,$$
which is in contrast with $\DTM(S) = \Mod_A$. 
Under this equivalence, the functor $M \mapsto M(1)$ (induced from $\DTM(S)$) corresponds to $N \mapsto N \twi {1}$, i.e., shifting the $\Z$-grading at the right hand category.

\item 
\label{item--dependence.on.S}
The notion of reduced motives depends on the choice of the base scheme $S$: for $X / S$, the natural map 
$$\DM(X) \t_{\DTM(S)} \grModL \r \DM(X) \t_{\DTM(\Spec \Z)} \grModL$$
usually won't be an equivalence.
Indeed, 
$$\DTM(S) \t_{\DTM(\Spec \Z)} \grModL \cong \Mod_{A_S} \t_{\Mod_{A_{\Spec \Z}}} \grModL = \Mod_{A_S \t_{A_{\Spec \Z}} \Lambda}$$
won't be equivalent to $\grModL$ since $A_S \t_{A_{\Spec \Z}} \Lambda \ne \Lambda$.

\item
\label{item--stacky.case}
For length reasons we focus our attention in this paper on reduced motives on ind-schemes.
For applications such as an integral Satake equivalence, one will need to construct categories such as $\MTMr(L^+G \backslash LG / L^+G)$, i.e.~mixed reduced Tate motives on the double quotient of the loop group by the positive loop group. 
This will require a category $\DMr(X)$ for any prestack $X$, or at the very least any algebraic stack $X$. 
Such a theory is available as soon as one has a robust (i.e., \ii-categorified) formalism of motives on (pre-)stacks. E.g., using the work in \cite{RicharzScholbach:Intersection}, a category $\DMr(X)_\Q$ with rational coefficients can be defined as in \refeq{DMr.definition}. 
Comparison results such as, say, \thref{comparison.graded} carry over unchanged to reduced motives in that generality.  
In order to capture genuine equivariant phenomena (with integral coefficients), one may use $\DM(X) := \Mod_{\Z_X}(\SH(X))$ for scalloped stacks, or for arbitrary stacks the limit-extended category $\SH_{\triangleleft}(X)$ defined by Khan and Ravi \cite{KhanRavi:Generalized}.
\end{enumerate}
\xrema

Recall that $\grModL$ is a symmetric monoidal \ii-category.
By definition, the symmetry isomorphisms obey the Koszul sign rule with respect to the cochain degree, while a shift in the graded direction (i.e., $\twi -$) causes no sign: the map
$$(\Lambda \twi n [m]) \t (\Lambda \twi {n'} [m']) \cong (\Lambda \twi {n'} [m']) \t (\Lambda \twi n [m])$$
is multiplication by $(-1)^{mm'}$. For example, the symmetric algebra over a $\Lambda$-module of the form $M \cong \Lambda^{\bigoplus I}[-1]\twi n$ is the \emph{exterior algebra} on $\Lambda^{\bigoplus I}$ (placed in degree $-1$).
The following lemma will allow us to do concrete computations with reduced Tate motives.

\lemm
\thlabel{DTM.flat.cell}
Let $X = \A_S^n \x_S (\GmX S)^m$.
Then there is an equivalence
$$\DTM(X) = \Mod_{\Sym (\Lambda(-1)[-1])^{\t m}} (\DTM(S))$$ 
and therefore an equivalence
$$\DTMr(X) = \Mod_{\Sym (\Lambda\twi {-1}[-1])^{\t m}}(\grModL).\eqlabel{DM.Gm.flat}$$
Thus, a reduced Tate motive on $X$ (with $\Lambda$-coefficients) can be colloquially described as a $\Z$-graded complex of $\Lambda$-modules $R$, together with $m$ anticommuting maps
$$R \r R\twi {+1}[+1] .$$
\xlemm

\pf
Let us abbreviate $E := (\Sym \Lambda (-1)[-1])^{\bigotimes m} = \Sym (\Lambda \twi {-1}[-1]^{\bigoplus m})$.
Let $q : X \r S$ be the structural map. The adjunction 
$$q^* : \DTM(S) = \Mod_A \rightleftarrows \DTM(X) : q_*$$ 
is monadic: $q_*$ is conservative since $q^* \Lambda(n)$ is a family of generators of $\DTM(X)$.
\nts{It exists, since $q_* \Lambda$ is a Tate motive. The right adjoint $q_*$ is conservative. Indeed, suppose $M \in \DTM(X)$ satisfies $q_* M = 0$.
Then $0 = \Map_{\DTM(S)}(\Lambda(n), q_* M) = \Map_{\DTM(X)}(q^* \Lambda(n), M) = 0$.
Since the $q^* \Lambda(n) = \Lambda_X(n)$ generate $\DTM(X)$, this shows $M = 0$. Moreover, $q_*$ preserves colimits. }
By the Barr--Beck monadicity theorem \cite[Theorem~4.7.3.5]{Lurie:HA}, $\DTM(X)$ is the category of algebras over the monad $\DTM(S)$ given by the endofunctor $q_* q^* - \stackrel \cong \gets q_* q^* \Lambda \t_\Lambda -$.
It remains to construct an isomorphism (in $\CAlg(\DTM(S))$)
$$\alpha : E \r q_* q^* \Lambda.$$
By the $\Sym$-forgetful adjunction, such a map is the same as a collection of $m$ maps (in $\DTM(S)$)
$$\Lambda(-1)[-1] \r q_* q^* \Lambda$$
or, yet equivalently, $m$ maps in $\DTM(X): \Lambda(-1)[-1] \r \Lambda$.
These maps arise via the *-projections along $\pr_i : X \r \GmX S$ from maps in $\DTM(\GmX S)$ of the form $\Lambda(-1)[-1] \r \Lambda$.
We may thus assume $m = 1$ to construct the map. 
The unit maps of adjunctions for *-pullbacks vs.~*-pushforwards, $\Lambda \r q_* q^* \Lambda \r p_* p^* \Lambda \cong \Lambda$ ($p : \A^1_S \r S$, cf.~\thref{DM.further.properties}) exhibit $\Lambda$ as a retract of $q_* q^* \Lambda$, with complement $\Lambda(-1)[-1]$. This yields a map $\alpha$ as stated.
It is then standard, using the Künneth formula for motivic cohomolgoy of $Z \x_S \GmX S$, that $\alpha$ is an isomorphism.

Under the equivalence $\DTM(S) = \Mod_A$, $E$ maps to $(\Sym A \twi {-1}[-1])^{\bigotimes m} = A \t_\Lambda (\Sym \Lambda \twi{-1}[-1])^{\bigotimes m} =: A \t F$.
Therefore, using generalities about tensor products of module categories \cite[Proposition~4.1]{BenZviFrancisNadler:Integral}
$$\Mod_{E}(\DTM(S)) = \Mod_{A \t F}(\Mod_A) = \Mod_{A} \t_{\grModL} \Mod_F(\grModL).$$
This implies the claim about $\DTMr(X)$.
\xpf
\rema \thlabel{DTM.flat.cell.koszul.dual} There is Koszul dual description of the category $\DTMr(X)$ for $X = \A_S^n \x_S (\GmX S)^m.$ Let $R := \Sym(\Lambda^m \twi 1)\in \grModL$ and consider the full subcategory $$\Modf_R(\grModL)\subset \Mod_R(\grModL)$$
generated by the objects $\Lambda\langle i \rangle$ for $i\in \Z$ by finite colimits and retracts.

Equivalently, $\Modf_R(\grModL)$ consists precisely of the $R$-modules whose underlying $\Lambda$-module lies in $\grPerf_\Lambda$. One direction of this statement is clear. To see the other direction, we claim that any complex $M$ of graded $R$-modules whose underlying graded $\Lambda$-module is perfect can be built from the simple $R$-modules $\Lambda(i)[j]$ inductively. For this, one starts with the highest non-zero degree part of $M$, say $N,$ which exists since $M$ is a perfect complex as a $\Lambda$-module. The action of $R$ on $M$ restricts to $N$ and factors through the augmentation map $R\to \Lambda$ for degree reasons. Hence, any resolution $N$ by the $\Lambda$-modules $\Lambda(i)[j]$ is also a resolution of $R$-modules, which proves the claim.

The objects $\Lambda(i)[i] \in \DTMr(X)^\comp$ generate a weight structure (which is distinct from the one in \thref{Chow.weight.structure}), whose heart is denoted by $\Cc$.
On the other hand, $\Cc$ is also the heart of a weight structure on $\Modf_R$ generated by the objects $\Lambda \twi i [i]$.  Since we have
$$\Hom_{\DMr(X)}(\Lambda,\Lambda(i)[i+n])=\Hom_{\Modf_R}(\Lambda,\Lambda\langle i\rangle[i+n])=0\text{ for all } n\neq 0$$
the weight complex functor provides equivalences of categories
$$\DTMr(X)^\comp\stackrel{\sim}{\to}\Comp^\bd(\Cc)\stackrel{\sim}{\leftarrow}\Modf_R(\grModL).$$
\xrema

The category of reduced motives has the following further immediate properties. In the entire lemma, $\DM$ can be replaced by $\DTM$ at will.

\lemm\phantom{a}
\thlabel{DM.flat.properties}
\begin{enumerate}[wide]
\item \label{item--DTM.flat.presentable}
$\DMr(X)$ is a presentable stable \ii-category. It is a module over $\grModL$, i.e., colloquially speaking, it is a $\Lambda$-linear category equipped with a $\Z$-grading.
In particular its homotopy category is a $\Z$-graded triangulated category having all coproducts and products.

\item
\label{item--flattening.functor}
There is a natural functor, called \emph{reduction functor}
$$\red : \DM(X) \r \DMr(X).\eqlabel{flattening.functor}$$

If $(M_i)_{i \in I}$ is a set of compact generators of $\DM(X)$ that is stable under applying Tate twists (in positive and negative direction), then the objects $\red(M_i)$ form a set of compact generators of $\DMr(X)$.
For two compact objects $M, M' \in \DM(X)$, we have
$$\IMaps_{\DMr(X)}(\red(M), \red(M')) = \IMaps_{\DM(X)}(M, M') \t_A \Lambda.\eqlabel{IMaps.DM.flat}$$
Here at the left $\IMaps$ denotes the enriched mapping object in $\grModL$ (by regarding $\DMr(X)$ as a $\grModL$-module), and the $\IMaps$ at the right denotes the enriched mapping object in $\Mod_A$, by virtue of $\DM(X)$ being a $\Mod_A$-module.
In particular,
$$\Maps_{\DMr(X)}(\red(M), \red(M')) = \ev_0 \left ( \IMaps_{\DM(X)}(M, M') \t_A \Lambda \right ).\eqlabel{Maps.DM.flat}$$

\end{enumerate}

\xlemm

\pf
\refit{DTM.flat.presentable}: This holds by the very definition of the Lurie tensor product.

\noindent\refit{flattening.functor}: The functor
$$\Lambda \t_A - : (\DTM(S) =) \Mod_A \r \grModL$$
is $\DTM(S)$-linear, where we regard the right hand category as a module over $\Mod_A$ via the augmentation map $A \r \Lambda$.
The reduction functor arises by applying $\DM(X) \t_{\DTM(S)} -$ to this functor.

The category $\Mod_A$ is a \emph{rigid} symmetric monoidal \ii-category with compact generators given by $A \twi n$ ($n \in \Z$). Thus, by \cite[Chapter~1, \S10.3, \S10.5.7]{GaitsgoryRozenblyum:StudyI}, objects of the form $M \boxtimes \Lambda \twi n \in \DMr(X)$ with $M$ running over a set of compact generators of $\DM(X)$ generate $\DMr(X)$.
The claim about the enriched mapping objects holds by Proposition~10.5.8 there.
\xpf

The description of generators and \refeq{Maps.DM.flat} immediately shows:

\coro
$\DTMr(X)$ is the presentable stable full subcategory of $\DMr(X)$ generated by $\red(\Lambda(n))$, $n \in \Z$.
Henceforth, we will denote these objects just by $\Lambda(n)$, if there is no ambiguity between $\DM(X)$ and $\DMr(X)$.
\xcoro

For later purposes we unwind the enriched mapping object $$\IMaps_{\DM(X)}(M, M') \in \Mod_A.$$
By definition, for each $r \in \Z$
$$\Map_{\Mod_A(\grModL)}(A \twi r, \IMaps_{\DM(X)}(M,M')) = \Maps_{\DM(X)}((A \twi r) \t M, M').$$
In the right hand side $A \twi r$ acts by a positive Tate twist, so the right hand side is $\Maps_{\DM(X)}(M, M'(-r))$.

The left hand side equals $\Map_{\grModL}(\Lambda \twi r, \IMaps_{\DM(X)}(M,M'))$, which equals the $(-r)$-th component of $\IMaps$.
In other words, 
$$\ev_r \IMaps_{\DM(X)}(M, M') = \Maps_{\DM(X)}(M, M'(r)).\eqlabel{IMaps.explicit}$$

\subsection{Functoriality for reduced motives}
\label{sect--DTMr.functoriality}

The formation $X \mapsto \DMr(X)$ is part of a six-functor formalism that we now describe.
In a nutshell, all functors ($f^*$, $f_*$, $f_!$, $f^!$, $\t$, $\Hom$) still exist for the categories $\DMr(X)$ and are compatible with the usual ones under the reduction functor.

Recall from \refsect{motives} that the formation $X \mapsto \DM(X)$ is part of a functor
$$\DM^*_! : \Corr := \Corr(\Sch_S^{\ft, \times, \dual})^\proper_{\sep, \all} \r \Mod_{\DTM(S)}.\eqlabel{DM.Corr}$$

\prop
\thlabel{flattening.functorial}\phantom{a}
\begin{enumerate}[wide] 
\item
\label{item--DM.flat.*!}
There is a lax symmetric monoidal functor
$$(\DMr)^*_! : \Corr \r \grPrStL$$
whose evaluation at $X / S$ is the category $\DMr(X)$ considered above.

\item
The reduction functors $\DM(X) \r \DMr(X)$ are then part of a natural transformation
$$\DM^*_! \r (\DMr)^*_!$$
between lax symmetric monoidal functors $\Corr(\Sch_S^\ft) \r \PrStL$; for concreteness we forget the $\Mod_A$-module structure on both functors at this point.
\end{enumerate}
\xprop

\rema
In parallel to \thref{Corr.explained}, the existence of the functor $(\DMr)^*_!$ encodes, in particular: for each map $f : X \r Y$ for $S$, there are functors 
$$f^* : \DMr(Y) \r \DMr(X),\eqlabel{f*.flat}$$
$$f_! : \DMr(X) \r \DMr(Y).\eqlabel{f!.flat}$$
Moreover, base-change and projection formulas as in \thref{Corr.explained} again hold for $\DMr$.
In addition, the second statement says that these functors are compatible with the usual *-pullback and !-pushforward under the reduction functors. 

In the same vein, reduced motives on ind-schemes are defined by composing the functor $\DM_!$ in \refeq{DM.indschemes} with $- \t_{\DTM(S)} \grModL$.
\xrema

\pf
We define $(\DMr)^*_!$ to be the composition
$$\Corr(\Sch_S^\ft) \stackrel{\DM^*_!} \lr \Mod_{\DTM(S)}(\PrSt) \stackrel{\grModL \t_{\DTM(S)} -} \lr \grPrStL (:= \Mod_{\grModL} (\PrSt)).$$
By definition, this reproduces the categories $\DMr(X)$ when evaluating at $X / S$.
Being the composite of the lax symmetric monoidal $\DM^*_!$ and the symmetric monoidal $\grModL \t_{\DTM(S)} -$, $(\DMr)^!_*$ is also lax symmetric monoidal.

The claimed natural transformation then results from the unit map of the adjunction
$$\grModL \t_{\DTM(S)} - : \Mod_{\DTM(S)} (\PrSt) \rightleftarrows \grPrStL$$
(and further applying the forgetful functor from the left hand category to $\PrStL$ induced by the map (of commutative algebra objects) $\ModL \stackrel{\twi 0} \r \grModL$).
\xpf
\prop
For each map of $S$-schemes $f : X \r Y$ the functors
$$f_* : \DM(X) \r \DM(Y),$$
$$f^! : \DM(Y) \r \DM(X)$$
are $\DTM(S)$-linear. The resulting functors
$$f_* \t_{\DTM(S)} \grModL, f^! \t_{\DTM(S)} \grModL$$
are the right adjoints of the functors in \refeq{f*.flat}, \refeq{f!.flat}.
\xprop

\pf
Since $\DTM(S)$ is rigid (being compactly generated by the dualizable objects $\Lambda(n)$, $n \in \Z$), this is an instance of \thref{adjunctions.tensor}.
\nts{One can work this out also slightly more hands on as follows: the functor $X \mapsto \DM(X)$, $f \mapsto f^*$ is part of a functor $\DM^* : (\Sch_S^\ft)^\opp \r \Alg_{LM}(\PrSt)$, which sends each $X$ to $\DM(X)$ regarded as a $\DTM(S)$-module.
We can conclude from \cite[Corollary~7.3.2.7]{Lurie:HA} the existence of an operadic right adjoint, since the evaluation of that functor in each object in $LM$ has a a right adjoint: on the ring-part the functor is just the identity on $\DTM(S)$, and on the module part it is the usual $f^*$, which does have $f_*$ as a right adjoint. 
In order to check that this right adjoint is not only a map of \ii-operads it suffices to check that the action maps such (for $f : X \r Y$, $x : X \r S$, $y : Y \r S$)
$$y^* A \t f_* B \r f_* (x^* A \t B)\eqlabel{f*.module.check}$$
are isomorphisms for all $B \in \DM(X)$ and $A \in \DTM(S)$. It suffices to check this for $A = 1(n)$, where it is standard.
In order to check that $f_* \t \id$ is a right adjoint of $f^* \t \id$ it suffices to observe that the unit and counit maps $\id \r f_* f^*$, $f^* f_* \r \id$ are $\DTM(S)$-linear. Then we use \thref{adjunctions.tensor} to conclude.}
\xpf

Given that the functors on reduced motives arise by base-changing the usual ones, the abstract features of $\DM$ carry over to $\DMr$.
Because of their importance for our purposes below, we specifically spell out homotopy invariance and localization.

\coro
(Homotopy invariance for reduced motives)
For a projection $p : \A^n_X \r X$, the functor
$$p^* : \DMr(X) \r \DMr(\A^n_X)$$
is fully faithful.
In particular, this restricts to an equivalence of categories
$$p^* : \DTMr(X) \stackrel \cong \r \DTMr(\A^n_X).$$
\xcoro

\pf
Indeed, the unit map $\id \r p_* p^*$ is an isomorphism since both $p_*$ and $p^*$ arise from their usual counterparts by base-changing along $\DTM(S) \r \grModL$.
\xpf

\coro
\thlabel{DM.flat.localization}
(Localization for reduced motives) 
For a diagram consisting of a closed immersion $i$ and its complementary open immersion $j$
$$Z \stackrel i \r X \stackrel j \gets U$$
we have a recollement situation, i.e., adjoints 
$$\eqalign{ j_! \dashv j^! & = j^* \dashv j_*, \cr
i^* \dashv i_* & = i_! \dashv i^!,}$$ 
such that $j_*$ and $i_*$ are fully faithful and such that $i^* j_! = i^! j_* = 0$.
\xcoro

\section{Stratified mixed Tate motives}
\label{sect--4}

In this section, we restrict the construction of reduced motives to specific geometric situations, namely stratified ind-schemes $X = \bigcup X_w$ where the $X_w$ are affine spaces or cells. We also restrict our attention to specific (reduced) motives, namely those motives $M$ such that all $M|_{X_w} $ are (reduced) Tate motives.

\subsection{Stratified Tate motives}

\label{sect--stratified.Tate}

In this section we recall some standard conventions on stratified (ind-)schemes, as in \cite[\S3]{RicharzScholbach:Intersection} or \cite[\S4]{SoergelWendt:Perverse} for schemes and define (reduced) stratified Tate motives.

\defi
A \emph{stratified ind-scheme} over $S$ is a map of ind-schemes over $S$
$$\iota\co X^+=\bigsqcup_{w\in W}X_w\,\r\, X$$
such that $\iota$ is bijective on the underlying sets, each stratum $X_w$ is a smooth $S$-scheme,
the restriction to each stratum $\iota_w:=\iota|_{X_w}:X_w\to X$ is representable by a quasi-compact immersion and
the topological closure of each stratum $\iota(X_w)$ is a union of strata. 

A stratification is \emph{cellular} (resp.~\emph{affine}) if each stratum $X_w$ is isomorphic to $\A^{n_w}_S \x \GmX S^{m_w}$ (resp., to $\A_S^{n_w}$).
\xdefi

\defilemm
\thlabel{DTM.stratified}
(\cite[\S4]{SoergelWendt:Perverse} for schemes, \cite[3.1.11]{RicharzScholbach:Intersection} for ind-schemes) 
A stratified (ind-)scheme is \emph{Whitney--Tate} if $\iota^* \iota_* \Lambda_{X^+} \in \DTM(X^+)$.
In this case, the following full subcategories of $\DM(X)^\comp$ are the same:
\begin{enumerate}
\item 
\label{item--upper.*}
The subcategory consisting of objects $M$ such that $\iota_w^* M \in \DTM(X_w)$ for all $w \in W$.

\item 
\label{item--upper.!}
The subcategory consisting of objects $M$ such that $\iota_w^! M \in \DTM(X_w)$ for all $w \in W$.

\item 
\label{item--lower.!}
The presentable stable subcategory generated by $(\iota_w)_! \Lambda(n)$ for all $w \in W$, $n \in \Z$.

\item 
\label{item--lower.*}
The presentable stable subcategory generated by $(\iota_w)_* \Lambda(n)$ for all $w \in W$, $n \in \Z$.
\end{enumerate}

We call this subcategory the category of \emph{stratified Tate motives}, denoted 
by $\DTM(X, X^+)$ or just $\DTM(X)$ if the presence of the stratification is clear from the context.
\xdefilemm

\defilemm \cite[\S3.1]{RicharzScholbach:Intersection}
Let $(X, X^+)$ and $(Y, Y^+)$ be two stratified (ind-)schemes and $f : X \r Y$ a schematic map of finite type, that is \emph{stratified} (i.e., $f$ is compatible with the stratification and maps each $X_w$ to some stratum $Y_{w'}$).
We say $f$ is a \emph{Whitney--Tate map}, if $f_*$ preserves stratified Tate motives ($f^*$ always does; cf.~\refsect{DM.ind-schemes} for the functor formalism of motives on ind-schemes).
Equivalently (since all strata are smooth over $S$), $f_!$ preserves Tate motives.
\xdefilemm

\exam \cite[Proposition 3.8]{EberhardtKelly:Mixed}
If a stratified map as above is such that for each $w$, the restriction $f|_{X_w}: X_w \r Y_{w'}$ is a surjective linear map between affine spaces, then $f$ is a Whitney--Tate map. 
We call such maps \emph{affine-stratified}.
\xexam

\defi
\thlabel{resolution.condition}
We say that an (ind-)scheme $(X, X^+)$ with an affine stratification \emph{admits affine-stratified resolutions} if for every stratum $X_w$ there is a resolution of singularities $\pi_w:\widetilde{X}_w\to \overline{X}_w$ (i.e., $\pi_w$ proper, $\widetilde X_w$ smooth over $S$) that is an isomorphism over $X_w$ and such that $\widetilde{X}_w$ admits an affine stratification and $\pi_w$ is affine-stratified. 
\xdefi

\rema 
The condition in \thref{resolution.condition} is motivated by similar conditions in \cite[Lemma 4.4.2]{Beilinson:Koszulduality}.
It ensures that an affine stratification is Whitney--Tate \cite[Proposition A.2]{SoergelWendt:Perverse}.
It will also play a rôle in considerations related to pointwise purity, see the proof of \thref{weight.complex.equivalence}.
\xrema
\exam The flag variety $X=G/B$ with its stratification in $B$-orbits fulfills the condition in \thref{resolution.condition}. The closure of $B$-orbits are Schubert varieties which have the Bott-Samelson resolution that is affine-stratified, see \cite{Haines:Purity}.
\xexam

\defilemm
\thlabel{DTM.flat.basic}
Let $\iota : X^+ \r X$ be a Whitney--Tate stratified (ind-)scheme. 
The category of \emph{reduced stratified Tate motives} is defined as 
$$\DTMr(X, X^+) := \DTM(X, X^+) \t_{\DTM(S)} \grModL.$$
If the choice $X^+$ is clear from the context, we abbreviate this by $\DTMr(X)$.

Equivalently, $\DTMr(X, X^+)$ is the full subcategory of $\DMr(X)$ characterized by the properties analogous to \thref{DTM.stratified}\refit{upper.*}--\refit{lower.*}, exchanging $\DTM$ by $\DTMr.$
The reduction functor $\DM(X)\to \DMr(X),$ see \refeq{flattening.functor}, restricts to a functor
$$\red : \DTM(X, X^+) \r \DTMr(X, X^+).$$
\xdefilemm

\pf
For an ind-scheme $X = \colim X_i$, we have $$\DTM(X, X^+) = \colim \DTM(\ol{X_i}, \ol{X_i}^+)$$ by \cite[Remark~3.1.3]{RicharzScholbach:Intersection}.
The transition functors in this filtered colimit are $\DTM(S)$-linear, so that we may assume $X$ is a scheme.

The description of generators of $\DTMr(X, X^+)$ as in \refit{lower.!} and \refit{lower.*} is a general consequence of the Lurie tensor product, as in \thref{DM.flat.properties}\refit{flattening.functor}.
The descriptions via the two pullback functors $\iota_w^*$ and $\iota_w^!$ is then a consequence of the localization formalism for reduced motives (\thref{DM.flat.localization}).
\xpf

\rema Denote by $\DMrx$ either $\DM$ or $\DMr.$
The subcategory of compact objects in $\DTMrx(X, X^+)$ is the subcategory of $\DMrx(X)$ generated, by means of finite colimits, shifts, and retracts by motives of the form $\iota_{w,*} \Lambda(n)$ (equivalently, $\iota_{w,!} \Lambda(n)$).

The category $\DTMrx(X, X^+)$ is compactly generated, i.e., $\DTMrx(X, X^+) = \Ind(\DTMrx(X,X^+)^{\comp})$, so that for many purposes it suffices to consider compact objects.
However, for an ind-scheme such as the affine Grassmannian $X = \Gr_G \stackrel p \r S$, the dualizing motive $\omega_{\Gr_G} := p^! \Lambda \in \DTMrx(\Gr_G)$ fails to be compact, so it is useful to have a presentable category of stratified Tate motives.
\xrema

Sheaves on stratified spaces can be concisely described as follows:

\rema
\thlabel{DM.lax.limit}
Sheaves on stratified spaces can be described using lax limits \cite[Example~4.1.6]{ArinkinGaitsgory:Category}: if $X^+ = U \sqcup Z \stackrel{j \sqcup i} \r X$ is a stratification by an open and a closed stratum, then
$$\DM(X) = \laxlim \left (\DM(U) \stackrel{i^* j_*} \r \DM(Z) \right ).$$
Indeed, this is an equivalent formulation of the localization property of motives (\thref{DM.further.properties}), according to which a motive on $X$ is equivalent to a triple $$(M_U, M_Z, M_Z \r i^* j_* M_U)$$ with $M_U \in \DM(U)$, $M_Z \in \DM(Z)$. 
The stratification is Whitney--Tate iff $i^* j_*$ preserves Tate motives, in which case we have
$$\DTM(X) = \laxlim \left (\DTM(U) \stackrel{i^* j_*} \r \DTM(Z) \right ).$$
We can then apply \thref{laxlim.tensor} and compute the category of stratified reduced Tate motives as
$$\DTMr(X) = \laxlim \left (\DTMr(U) \stackrel{i^* j_*} \r \DTMr(Z) \right ).$$
For example, objects in $\DTMr(\P, \A^1 \sqcup \{\infty\})$ are triples 
$$\left (M, M' \in \grModL, M' \r M \t (\Lambda \oplus \Lambda \twi{-1}[-1])\right ).$$
\xrema

\subsection{Perverse t-structures}
In this section, we define a perverse t-structure on the categories of (reduced) stratified Tate motives on (ind-)schemes with cellular Whitney--Tate stratification.
We begin with the case of a cell itself. The definition follows the usual convention for perverse sheaves on complex varieties, that is, it makes use of the middle perversity.

\defilemm \thlabel{DTM.flat.tstructureonstrata}
Let $X=\A^n_S \x_S \GmX S^m.$ There is a unique t-structure, referred to as the \emph{perverse $t$-structure}, on $\DTMr(X)$ with heart 
$$\MTMr(X) := \DTMr(X)^{\tt = 0} \subset \DTMr(X)$$ 
generated by coproducts and extensions of the objects $\Lambda (q) [n+m]$ for $q\in \Z.$
This t-structure restricts to a t-structure on the subcategory of compact objects.

The statement holds true for $\DTM(X)$ if 
\begin{enumerate}
  \item \label{item--Lambda.BS} $\Lambda=\Q$ and
  \item \label{item--S.BS} $S$ is the spectrum of a finite field, or a global field or the ring of integers in a global field.
\end{enumerate}
\xdefilemm

\pf
For $\DTM(X)_\Q$ the statement follows from \cite[Theorem 4.2]{Levine:Tate}. The Beilinson--Soul\'e vanishing condition on  $\Hom_\DM(X)(\Q, \Q(i)[n])$ follows from the computations of algebraic K-theory of $S$ due to Quillen, Borel, and Harder. 

For $\DTMr(X),$ one can proceed along similar lines, using the computation of $\Hom_{\DTMr(X)}(\Lambda, \Lambda(i)[n])$ in \thref{DTM.flat.cell}.
Alternatively, one may deduce the statement using the equivalence between $\DTMr(X)^\comp$ and the category $\Modf_R(\grModL)$ for $R := \Sym(\Lambda^m \twi 1)$ mapping  $\Lambda(i)$ to $\Lambda\twi i,$ see \thref{DTM.flat.cell.koszul.dual}.
Recall that the category $\Modf_R(\grModL)$ consists of objects whose underlying $\Lambda$-module is perfect. Hence, it inherits a $t$-structure such that the forgetful functor to $\grPerf_\Lambda$ is $t$-exact (with respect to the usual t-structure on $\grPerf_\Lambda$, cf.~\thref{ModL.t.weight.structure}).
This yields the existence of the desired $t$-structure on $\DTMr(X)^\comp.$ This is the restriction of $t$-structure on $\DTMr(X)$ using that this category is compactly generated by the objects $\Lambda(i)$ in the heart.
\xpf

\rema 
The $t$-structure on $\DTM(X)$ also restricts to compact objects for more general coefficient rings $\Lambda.$ This is work in progess of Spitzweck--Uschogov, see also \cite[Theorem 9.10.]{Spitzweck:Mixed}.
\xrema

\defilemm
\thlabel{DTM.flat.t-structure}
Let $(X, X^+)$ be an (ind-)scheme with a cellular Whitney--Tate stratification.
The \emph{perverse $t$-structure} on $\DTMr(X, X^+)$ is the $t$-structure glued from the perverse t-structures on the categories $\DTMr(X_w)$:
\begin{align*}
  \DTMr(X,X^+)^{\tt\leq 0}&=\{ M\in \DTMr(X,X^+)\mid \iota^* M\in \DTMr(X^+)^{\tt\leq 0}\}\text{ and }\\
  \DTMr(X,X^+)^{\tt\geq 0}&=\{ M\in \DTMr(X,X^+) \mid \iota^! M\in \DTMr(X^+)^{\tt\geq 0}\}.
\end{align*}
This t-structure restricts to a t-structure on the subcategory $\DTMr(X, X^+)^\comp$ of compact objects.
Again, under the conditions \ref{DTM.flat.tstructureonstrata}\refit{Lambda.BS} and \refit{S.BS} above, the same statement holds for $\DTM(X)$.
In this event, the reduction functor $\red : \DTM(X) \r \DTMr(X)$ is t-exact. We denote the respective hearts of this t-structure by $$\MTMrx(X,X^+)\subset\DTMrx(X, X^+)$$ and refer to it as the category of \emph{stratified mixed (reduced) Tate motives}. 
\xdefilemm

\pf
The t-structure is a routine consequence of the gluing formalism from \cite{BeilinsonBernsteinDeligne:Faisceaux}.
The t-exactness of $\red$ holds since $\red$ commutes with $\iota^*$ and $\iota^!$.
\xpf





In order to address whether the realization functor \refeq{realization.functor}
\[\D^\bd(\MTMr(X,X^+)^\comp)\to \DTMr(X,X^+)^\comp\]
is an equivalence, we use tilting objects.
This formalism was developed in the context of highest weight categories in \cite{Ringel:CategoryModulesGood1991}. See \cite{Beilinson:TiltingExercises2004} for an account of tilting objects in the geometric setting. To apply this theory, it is necessary to have standard and costandard objects in the category $\MTMr(X,X^+)^\comp,$ which necessitates that the functor $\iota_*$ is $t$-exact. In the context of $\ell$-adic sheaves this is true by Artin vanishing, which implies that pushforward along affine maps is $t$-exact. In our context we have to impose this as an additional condition. 

\defi 
\thlabel{condition.dagger}
Let $(X,X^+)$ be an (ind-)scheme with an affine Whitney--Tate stratification such that the functor $\iota_*$ is $t$-exact.
Let $\Lambda$ be a PID (principal ideal domain, e.g., a field).
\begin{enumerate}
  \item 
  We call the objects 
  \[\Delta_w(i)= i_{w,!}\Lambda(i)[\dim X_w]\text{ and }\nabla_w(i)= i_{w,*}\Lambda(i)[\dim X_w]\] the \emph{standard} and \emph{costandard objects}, respectively.
  \item An object $T\in \MTMr(X,X^+)$ is called a \emph{tilting object} if it admits finite filtrations such that the associated graded objects are finite direct sums of direct summands of standard and costandard objects, respectively. Tilting objects span a full subcategory denoted by $\Tilt(X) \subset \MTMr(X, X^+)^\comp$.
\end{enumerate}
\xdefi

\lemm 
\thlabel{DTM.noextbetweenstdandcostandard} 
Let $(X,X^+)$ be as in \thref{condition.dagger}.
If $M,N\in \MTMr(X)^\comp$ admit a standard and costandard filtration, respectively, then
  \[\Hom_{\DMr(X)}(M,N[j])=0 \text{ for all } j\neq 0.\]
\xlemm

\pf 
The group
$\Hom_{\DMr(X)}(\iota_{v,!}\Lambda[dim(X_v)],\iota_{w,*}\Lambda[\dim(X_w)](i)[j])$ can be computed as $\Hom_{\DMr(X_v)}(\Lambda,\iota_{v}^!\iota_{w,*}\Lambda[\dim(X_w)-\dim(X_v)]](i)[j])$.
This vanishes if $v \ne w$ by base change. For $v = w$ it vanishes for $j \ne 0$ since $\DTMr(X_v) = \grModL$.
The statement follows by an induction on the length of the filtration. 
\xpf

\prop\thlabel{DTM.enoughtiltingperverse} 
Let $(X,X^+)$ be as in \thref{condition.dagger} and let $\Lambda$ is a PID.
\begin{enumerate}
\item 
\label{item--tilting.exists}
For each stratum $\iota_w:X_w\to X,$ there is a tilting object $\Tt_w\in \Tilt(X)$ supported on $\overline{X}_w$ and $\iota_w^*\Tt_w=\Lambda[dim(X_w)]$.
\item 
\label{item--tilting.equivalence}
The realization functor \refeq{realization.functor} is an equivalence of categories
  \[\D^\bd(\MTMr(X,X^+)^\comp) \stackrel \cong \to \DTMr(X,X^+)^\comp.\]
\end{enumerate} 
\xprop
\pf 
\refit{tilting.exists}: We may assume $X$ is a scheme since the closure $\overline{X}_w$ of any stratum is a scheme and the functor $\ol \iota_{w,!}:\DTMr(\overline{X}_w) \to \DTMr(X)$ for $\ol \iota_w:\overline{X}_w\to X$ is $\tt$-exact and preserves tilting objects. Now, the argument in \cite[Proposition B.3]{AcharRiche:ModularI} for constructible sheaves on complex varieties, stratified by affine spaces, with integral coefficients translates unchanged to our setting. 

\refit{tilting.equivalence}:
The realization functor is fully faithful on the additive subcategory $\Cc$ of $\D^\bd(\MTMr(X)^\comp)$ generated by 
$\bigcup_{n\in\Z}\Tilt(X)[n]$: 
by \thref{DTM.noextbetweenstdandcostandard} there are no non-trivial extensions between tilting objects in $\DTMr(X)$ and thus in particular there are no non-trivial extensions in $\MTMr(X).$
By the five lemma the realisation functor is also fully faithful on the stable category generated by $\Cc$. By \refit{tilting.exists}, tilting objects generated $\DTMr(X)^\comp$.
Therefore they also generate $\MTMr(X)^\comp$, and hence $\D^\bd(\MTMr(X)^\comp)$.
\xpf

\rema  Weight structures provide a convenient perspective on tilting objects.
 Namely, there is a weight structure $\wt'$ called \emph{tilting weight structure} on $\DTMr(X,X^+)^{\comp}$ such that $$\Tilt(X)=\MTMr(X,X^+)^\comp \cap \DTMr(X,X^+)^{\comp, \wt'=0}.$$
For a single stratum $X=\A^n_S$, the heart of the tilting weight structure on $\DTMr(X)^\comp$ is generated by the objects $\Lambda(r)[n]$ for $r\in\Z$ under finite direct sums and direct summands.
In general, the tilting weight structure on $\DTMr(X,X^+)^\comp$ is obtained by gluing. 
We refer to \cite[Section 2]{Achar:Cotstructures} for a very thorough discussion on the relation of perverse $t$-structure and the tilting weight structure.

In the situation of \thref{DTM.enoughtiltingperverse}, the realization functor for the perverse $\tt$-structure and the weight complex functor for tilting weight structure $\wt'$ yield a chain of equivalences
\[\D^\bd(\MTMr(X,X^+)^\comp) \stackrel \cong \to \DTMr(X,X^+)^\comp \stackrel \cong \to \Comp^\bd(\Tilt(X)).\]
\xrema

\subsection{Weight structures}
Similarly to mixed $\ell$-adic sheaves, there is a notion of weights on the category of motives $\DM(X).$ Weight structures yield a convenient language for this yoga of weights. 
The \emph{Chow weight structure} on $\DM(X)$ is the weight structure whose heart is generated by direct summands of motives $\pi_*\Lambda$ where $Y$ is regular and $\pi:Y\to X$ is proper. 
This necessitates resolution of singularities, and thus entails restrictions on the coefficient ring $\Lambda$, e.g. for $S = \Spec \Z$ one needs to take $\Lambda = \Q$.
For schemes with cellular stratifications, such as the ones appearing in geometric representation theory, one can construct weight structures for all $\Lambda$ more directly.
As before, let $\DTMrx$ denote either the category of Tate motives or the category of reduced Tate motives.

\defilemm 
\thlabel{Chow.weight.structure}
Let $X=\A^n_S \x_S \GmX S^m$ or a disjoint union of such schemes. There is a unique weight structure, referred to as the \emph{Chow weight structure}, on $\DTMrx(X)$ whose heart 
$$\DTMrx(X)^{\wt = 0} \subset \DTMrx(X)$$ 
is generated by direct sums and direct summand by the objects $\Lambda (q) [2q]$ for $q\in \Z.$
This weight structure is Ind-extended (\thref{weights.Ind}) from a weight structure on the subcategory of compact objects.
For example the heart $\DTMr(X)^{\wt = 0}$ of reduced motives on $X=\A^n_S$ is equivalent to the ordinary category of graded projective $\Lambda$-modules.
\xdefilemm

\pf
To see this, one has to show that 
$\Hom_{\DTMr}(X)(\Lambda,\Lambda(q)[2q][i])=0$ for all $q\in\Z$ and $i>0.$ This follows from the computation of motivic cohomology of $X$, as in \thref{DTM.flat.cell}.
\xpf

As for the perverse $t$-structure in \thref{DTM.flat.t-structure} the Chow weight structure on $\DTMrx(X, X^+)$ can be obtained via gluing and the reduction functor preserves weights since it commutes with $\iota^*$ and $\iota^!$.

\defilemm
\thlabel{DTM.flat.weightstructure}
Let $(X, X^+)$ be a cellular Whitney--Tate stratified (ind-)scheme.
The \emph{Chow weight structure} on $\DTMrx(X) := \DTMrx(X, X^+)$ is the weight structure glued from the Chow weight structures on the strata $\DTMrx(X_w).$ That is
\begin{align*}
  \DTMrx(X,X^+)^{\wt\leq 0}&=\{ M\in \DTMrx(X,X^+)\mid \iota^* M\in \DTMrx(X^+)^{\wt\leq 0}\}\text{ and }\\
  \DTMrx(X,X^+)^{\wt\geq 0}&=\{ M\in \DTMrx(X,X^+)\mid \iota^! M\in \DTMrx(X^+)^{\wt\geq 0}\}.
\end{align*}
This weight structure is Ind-extended from a weight structure on the subcategory $\DTMrx(X,X^+)^\comp$ of compact objects. The reduction functor $\red : \DTM(X,X^+) \r \DTMr(X,X^+)$ is weight exact.
\xdefilemm
Next, we compare the heart of the Chow weight structure for reduced and non-reduced Tate motives on a point.
\prop
\thlabel{heart.equivalence}
The following are equivalent:
\begin{enumerate}
\item 
\label{item--no.higher.Chow}
The Chow groups $\CH^n(S, \Lambda)$ with $\Lambda$-coefficients vanish for $n > 0$.
(This is the case if for example $S = \Spec k$ for a field $k$ or $S = \Spec \Z$.)

\item
\label{item--flat.eq.hearts}
The restriction of the reduction functor to weight-zero objects,
$$\red: \Ho(\DTM(S)^{\wt = 0}) \r \Ho(\DTMr(S)^{\wt=0})$$
is an equivalence of (additive) categories.
\item
\label{item--flat.ff}
For any $M^{\le 0} \in \DTM(S)^{\wt \le 0}$ and $M^{\ge 0} \in \DTM(S)^{\wt \ge 0}$, the map 
$$\Hom_{\DM(S)}(M^{\le 0}, M^{\ge 0}) \r \Hom_{\DMr(S)}(\red(M^{\le 0}), \red(M^{\ge 0}))$$
is an isomorphism.
\end{enumerate}

\xprop

\pf
In \refit{flat.eq.hearts}, being an equivalence is equivalent to being fully faithful since the generators $\Lambda(n)[2n]$ are in the image by design.
This immediately shows \refit{flat.ff} $\Rightarrow$ \refit{flat.eq.hearts}.

Condition \refit{flat.ff} is equivalent to having an isomorphism
$$\Hom_{\DM(S)}(\Lambda, \Lambda(n)[2n+i]) \stackrel \cong \r \Hom_{\DMr(S)}(\Lambda, \Lambda(n)[2n+i])\eqlabel{Hom.claim}$$
for all $n \in \Z$ and all $i \ge 0$.
Indeed, by \thref{weights.Ind}\refit{Ind.C.weights.description}, $\DTM(S)^{\wt \le 0}$ is the smallest subcategory stable under extensions and coproducts and containing $\Lambda(n)[2n+i]$ for $i \le 0$.
The dual description for $\DTM(S)^{\wt \ge 0}$ and the compactness of $\Lambda \in \DTM(S)$ reduces us to this special case.

The right hand side in \refeq{Hom.claim} identifies with $\Hom_{\grModL}(\Lambda, \Lambda \twi {n} [2n+i])$, which vanishes for all $n \ne 0$ (even for all $i \in \Z$).
The left hand side always vanishes for $i > 0$, cf.~\thref{DM.further.properties}\refit{algebraic.cycles}.
This shows all the remaining implications \refit{no.higher.Chow} $\Leftrightarrow$ \refit{flat.eq.hearts} $\Rightarrow$ \refit{flat.ff}, since $\Hom_{\DM(S)}(\Lambda, \Lambda(n)[2n]) = \CH^n(S, \Lambda)$ always vanishes for $n < 0$.
\xpf
In the following discussion we will need to impose an additional pointwise purity condition.
\defi For $?\in\{*,!\}$, a motive $M\in \DTMrx(X, X^+)$ is called $?$-\emph{pointwise pure} if $\iota^{?}M\in \DTMrx(X^+)^{\wt =0}$.
\xdefi
\prop \thlabel{heart.noexts}
Let $(X,X^+)$ be an (ind)-scheme with an affine stratification. Let $M,N\in \DTMr(X,X^+)^{\wt=0}$ be $*$-pointwise and $!$-pointwise pure, respectively. Then for all $i\neq 0$
$$\Hom_{\DMr(X)}(M,N[i])=0.$$
\xprop
\pf The case of a single stratum $X=\A^n_S$ follows from the explicit description $\DTMr(\A^n_S)\cong\DTMr(S)=\grMod_\Lambda$. The general case follows by induction as in \cite[Lemma~3.16]{EberhardtKelly:Mixed}.
\xpf

\prop \thlabel{heart.fullyfaithful}
Assume that $\CH^n(S, \Lambda)=0$ for $n>0$. Let $(X,X^+)$ be an (ind)-scheme with an affine stratification and $M,N\in \DTMr(X,X^+)^{\comp,\wt=0}$ be $*$-pointwise and $!$-pointwise pure, respectively. Then reduction gives an isomorphism
$$\Hom_{\DM(X)}(M,N)\to \Hom_{\DMr(X)}(\red(M),\red(N)).$$
\xprop

\pf The case of a single stratum $X=\A^n_S$ follows from \thref{heart.equivalence} since $\DTMrx(\A^n_S)\cong\DTMrx(S)$ by homotopy invariance.
We may replace $X$ by the support of $M$ and $N$ which is a scheme given that they are compact object.
Let
$i:Z\to X\gets U:j$ be the inclusion of a closed stratum and its open complement. Then, the localisation sequence yields the following diagram of exact sequences of $\Hom$-groups:
\begin{center}
\begin{tikzcd}[column sep=5pt]
\Hom(j^*M,j^!N[-1]) \arrow[d] \arrow[r] & \Hom(i^*M,i^!N) \arrow[d] \arrow[r] & \Hom(M,N) \arrow[d] \arrow[r] & \Hom(j^*M,j^!N) \arrow[d] \arrow[r] & 0 \arrow[d] \\
0 \arrow[r]     & \Hom_{\red}(i^*M,i^!N)\arrow[r]     & \Hom_{\red}(M,N)\arrow[r]     & \Hom_{\red}(j^*M,j^!N) \arrow[r]     & 0.
\end{tikzcd}
\end{center}
Here we abbreviated $\Hom_{\red}(X,Y)=\Hom_{\DMr}(\red(X),\red(Y)).$ The two zeroes in the right column come from the axioms of a weight structure. The zero in the bottom left follows from \thref{heart.noexts}. The second vertical arrow is an isomorphism by \thref{heart.equivalence}. The fourth vertical arrow is an isomorphism by induction using that $?$-restriction preserves $?$-pointwise purity. The five lemma implies that the third vertical arrow is also an isomorphism.
\xpf

\prop 
\thlabel{weight.complex.equivalence}
Assume that $\CH^n(S, \Lambda)=0$ for $n>0$. 
Let $(X,X^+)$  be an (ind-)scheme with an affine stratification that admits affine-stratified resolutions, see \thref{resolution.condition}.
Then the reduction functor $\red$ is equivalent to the weight complex functor on $\DTM(X)$. More precisely, we have a commutative diagram with equivalences as indicated:
\begin{center}
  \begin{tikzcd}
    {\DTM(X,X^+)^\comp} \arrow[r] \arrow[d,"\red"] & {\Ch^\bd(\Ho(\DTM(X,X^+)^{\comp,\wt=0}))} \arrow[d,"\wr"] \\
    {\DTMr(X,X^+)^{\comp}} \arrow[r,"\sim"]                            & 
    {\Ch^\bd(\Ho(\DTMr(X,X^+)^{\comp,\wt=0})).}                                          
  \end{tikzcd}
  \end{center}
\xprop
\pf 
The diagram is commutative by the functoriality of the weight complex functor \cite[Corollary~3.0.4]{Sosnilo:Theorem}.

Choose affine-stratified resolutions $\pi_w:\widetilde{X}_w\to\overline{X}_w.$ Let $\Ee_w=\pi_{w,!}\Lambda.$ Then $\Ee_w$ is supported on $\overline{X}_w.$ By an argument as in \cite[Theorem 4.5]{EberhardtKelly:Mixed} the objects $\Ee_w(n)[2n]$ for $n\in \Z$ are pointwise pure and generate $\DTMrx(X,X^+)^\comp$ under finite colimits and retracts and generate $\DTMrx(X,X^+)^{\comp,\wt=0}$ under finite direct sums and retracts. Now \thref{heart.fullyfaithful} implies that the right arrow is an equivalence. The bottom arrow is an equivalence using \thref{heart.noexts}, the fact that both sides are generated by the $\Ee_w(n)[2n]$ and the five lemma.\xpf

\subsection{Independence of the base scheme}

A useful feature of reduced motives being defined in arbitrary characteristic is that it becomes possible to switch the ground scheme $S$, for example one can mediate between characteristic 0 and characteristic $p$ base schemes.

\prop
\thlabel{independence.S}
Consider a cartesian diagram
$$\xymatrix{
X'^+ \ar[r]^{\iota'} \ar[d]^{s^+} & X' \ar[d]^s \ar[r] & S' \ar[d]^{s_0} \\
X^+ \ar[r]^\iota & X \ar[r] & S.
}$$
Here $s_0$ is a map of base schemes and $\iota$ is a cellular Whitney--Tate stratification on an (ind-)scheme $X$, and $\iota'$ gives the pulled back stratification on $X' := X \x_S S'$.
Assume that the stratifications are such that the natural map
$$s^* \iota_* \r \iota'_* (s^+)^*\eqlabel{exchange.stratification}$$
is an isomorphism of functors.

Then the stratification $\iota'$ on $X'$ is again cellular Whitney--Tate and the natural functor 
$$s^* : \DTMr(X) \r \DTMr(X')$$
is an equivalence.
Here the reductions refer to the respective base schemes, i.e. over $S$ for $X$ and over $S'$ for $X'$ (cf.~\thref{DTM.flat.examples}\refit{dependence.on.S}).
\xprop

\exam 
The condition \refeq{exchange.stratification} holds for partial (affine) flag varieties \cite[Lemma~2.12]{RicharzScholbach:Motivic}. Moreover, it can be shown that 
the conclusion ot the theorem hold for all (ind-)schemes with an affine stratification that admit affine-stratified resolutions, see \thref{resolution.condition}.
\xexam

\rema
This statement can be compared with, say, the independence of the motivic Satake category of the base scheme \cite[Corollary~6.6.]{RicharzScholbach:Motivic}.
In both cases a certain semi-simplification has been performed, for $\DMr(X)$ in the guise of applying $- \t_A \Lambda$; for the Satake category by definition.
The difference is that the construction here works for stable \ii-categories of motives, which is useful for applications involving the full stable (or triangulated) category of sheaves on, say, $G/B$, as opposed to the subcategory of perverse sheaves. 
\xrema

\pf
The category $\DTM(X')$ is generated by objects of the form $(\iota'_w)_! 1(n)$, for a stratum $\iota'_w : X'_w \r X$.
These are clearly in the image of $s^*$, so it remains to check full faithfulness.

The Whitney--Tate condition for $X'$, i.e., $\iota'_* \iota'^* \Lambda_{X'} \in \DTM(X', X'^+)$ holds since $\iota'_* \iota'^* \Lambda = \iota'_* \iota'^* s^* \Lambda \stackrel{\refeq{exchange.stratification}} = (s^+)^* \iota_* \iota^* \Lambda$, which is in $\DTM(X'^+)$ since $\iota_* \iota^* \Lambda \in \DTM(X^+)$.

As for the claimed equivalence, first suppose $X$ is a cell, i.e., $X = \A^n_S \x (\GmX S)^m$. In this case the claim holds by \thref{DTM.flat.cell}.
Then, an induction reduces us to the case of two strata $j : X_0 \subset X$ (open) and $i : X_1 \subset X$ (closed), in which case we use \thref{laxlim.tensor} (cf.~\thref{DM.lax.limit}):
$$\DTMr(X, X^+) = \laxlim \left (\DTMr(X_0) \stackrel{i^* j_*} \lr \DTMr(X_1) \right).$$
The functor $s^*$ induces an equivalence on each term in the lax limit, and moreover commutes with $i^* j_*$ by assumption.
Thus it induces an equivalence $\DTM(X, X^+) \r \DTM(X', X'^+)$ as claimed.
\xpf

\section{Comparison results}
\label{sect--comparison.results}
In this section, we show that the category of reduced stratifies Tate motives recovers, refines and unifies the existing approaches to mixed sheaves in the literature, namely (unreduced) Tate motives over $S = \Spec \Fq$ (for $\Lambda = \Q$ or $\Fp$), semisimplified Hodge motives (over $S = \Spec \C$), graded $\ell$-adic sheaves (over $S = \Spec \Fp$) and Achar--Riche's mixed category (over $S = \Spec \C$).
In addition to the discussion after \thref{DTM.comparisons}, let us include a few more technical comments: the usage of $\Qlq$-adic sheaves in Ho--Li's approach of course allows to use the whole $\ell$-adic arsenal including the existence of perverse t-structures and Artin vanishing.
The motivic approaches in \cite{EberhardtKelly:Mixed,SoergelWendt:Perverse}, as well as the one presented here, typically require a closer look at the geometric objects at hand.
The restriction to (ind-)schemes with a cellular or affine stratification or presence of the condition that the pushforward along inclusions of strata is $t$-exact (see \thref{condition.dagger}) are a consequence of this state of affairs.
However, for applications in geometric representation theory this is no drawback since the combinatorics encoded in the geometry of partial (affine) flag varieties eventually requires using such properties anyways.

Our comparison results allow to bridge the gap between these different approaches in the literature:

\coro
Let $X / \Spec \Z$ be a cellular stratified Whitney--Tate (ind-)scheme with special fiber $X_p := X \x_{\Spec \Z} \Spec \Fp$ and generic fiber $X_\C := X \x_{\Spec \Z} \Spec \C$. 
Suppose that the condition in \refeq{exchange.stratification} is satisfied for $s : \Spec \Fp \r \Spec \Z$.
Then a choice of an isomorphism $\Qlq \cong \C$ yields equivalences
$$\Shv_{\Tate,\gr}(X_p) \cong \DTM(X_p,X^+_p)_\C \cong \DTM_\calH(X_\C,X^+_\C)\cong \DTMr(X,X^+)_\C$$
If the stratification is affine and admits affine-stratified resolutions, then there are equivalences
$$\DAR(X_\C)_\Fp \cong  \DTM(X_p,X^+_p)^\comp_\Fp\cong \DTMr(X,X^+)^\comp_\Fp$$
\xcoro

\exam
The assumptions on $X$ and $s$ are satisfied for $X$ a partial flag variety $G / P$ or a partial affine flag variety such as the affine Grassmannian $\Gr_G$ or the affine flag variety $\Fl_G$ associated to a reductive group $G / \Spec \Z$.
\xexam

\pf
The condition in \refeq{exchange.stratification} is automatic for $\eta : \Spec \C \r \Spec \Z$, since it is a pro-étale map.
Thus, \thref{comparison.Hodge}, \thref{comparison.graded}, and \thref{independence.S} yield a number of equivalences of categories
$$\xymatrix{
\DTMr(X_p)_\Qlq \ar[d]^\cong & \DTMr(X)_\Qlq \ar[l]^{s^*}_\cong \ar[r]^\cong & \DTMr(X)_\C \ar[r]_{\eta^*}^\cong & \DTMr(X_\C)_\C \ar[d]^\cong \\
\Shv_{\Tate, \gr}(X_p, \Qlq) & & & \DTM_\calH(X_\C). 
}$$
The second chain of equivalences similarly follows from \thref{comparison.DM} and \thref{DTM.flat.paritysheaves}
\xpf

\subsection{Motives over finite fields}

\prop
\thlabel{comparison.DM}
Let $S = \Spec \Fq$ or $\Spec \Fqq$ and $\Lambda = \Q$ or $\Lambda = \mathbf{F}_{p}$.
For any scheme $X / S$, the reduction functor is an equivalence of categories:
$$\red : \DM(X)_\Lambda \stackrel \cong \r \DMr(X)_\Lambda.$$
The same holds true for stratified Tate motives on stratified ind-schemes.
\xprop

\pf
For $\Lambda = \Q$ or $\Fp$ the unit map $\Lambda \r A$ is a (graded) quasi-isomorphism, i.e., $\H^n(A_r) = \H^n(S, \Lambda(r)) = 0$ for $r \ne 0$ or $n \ne 0$.
Indeed, for $\Lambda = \Q$, this group is isomorphic to $K_{2r-n}(S)^{(r)}_\Q$ which vanishes by Quillen's computation of K-theory of finite fields (and continuity of K-theory in case $S = \Spec \Fqq$).
For $\Lambda = \Fp$, this again vanishes as a consequence of Geisser--Levine's computation of mod-$p$ motivic cohomology, see \cite[Corollary~2.53]{EberhardtKelly:Mixed}.
By the above quasi-isomorphism the reduction functor $\DTM(S, \Lambda) \r \grMod_\Lambda$ is an equivalence, giving our claim.
\xpf

\rema
For any field $k \ne \Fq, \Fqq$, the reduction functor $\red : \DM(\Spec k)_\Q \r \DMr(\Spec k)_\Q$ is not an equivalence, since $\H^1(k, \Q(1)) = k^\x \t \Q \ne 0$.
\xrema

\subsection{Comparison with semisimplified Hodge motives}

Let $S = \Spec \C$ and $\Lambda=\C$.
We will show that reduced Tate motives with complex coefficients reproduce the category of \emph{semisimplified Hodge motives} due to Soergel and Wendt \cite{SoergelWendt:Perverse,SoergelVirkWendt:Equivariant}.

Recall the functor
$$\SmAff^\opp / S \r \Comp(\Ind(\MHS^\pol_\Q)) \stackrel{\gr^\W} \r \Comp(\Ind(\HS^{\pol, \Z}_\C)).$$
The first functor maps any $X / S$ to a complex of mixed Hodge structures whose $n$-th cohomologies are Deligne's mixed Hodge structures on the Betti cohomology of the associated complex manifold $X^\an$, $\H^n(X^\an, \Q)$ \cite{Drew:Realisations}. 
The functor $\gr^\W : \MHS^\pol_\Q \r \gr\HS^{\pol}_\C$ takes a (polarized) mixed Hodge structure and associates to it the graded pieces of the weight filtration \cite[Proposition~2.11]{SoergelWendt:Perverse}. This is an exact $\t$-functor (this uses complex coefficients).
By op.~cit., this composite functor passes to a symmetric monoidal colimit-preserving functor, called the \emph{semisimplified Hodge realization functor}
$$\R_\calH: \DM(S) \r \D(\Ind(\gr\HS^{\pol}_\C)).$$
This functor has a right adjoint $\R_*$ which one uses to define
$$\calH := \R_* \C \in \DM(S), \calH_X := f^* \calH \in \DM(X)_\C$$
for any scheme $f : X \r S$.
By the setup, $\calH_X$ is a commutative algebra object in $\DM(X)$, and one can consider $\DM_\calH(X) := \Mod_{\calH}(\DM(X))$.
This category is called the category of \emph{semisimplified Hodge motives}. 
Just as for $\DM$, there is a six-functor formalism for $\DM_\calH$ and a concomitant category of (stratified) Tate motives, which we denote by $\DTM_\calH \subset \DM_\calH$.

\prop
\thlabel{reduced.Hodge.realization}
Let $X / S$ be a scheme. Then the functor $\R_\calH$ induces a functor
$$\R_{\calH, \red}: \DMr(X) \to \DM_\calH(X)$$
compatible with the six functors on both sides.
\xprop

\pf
Let $A$ be defined as in \refeq{A}. We will use the following remark: since the unit map $1 : \C \r A$ is a quasi-isomorphism in graded degree 0
any morphism $A \r \C$ (in $\CAlg(\grModL)$) is uniquely determined by its restriction along the unit map, where it corresponds to a ring homomorphism $\C \r \C$ within the ordinary category of $\C$-vector spaces.

According to \cite[p.~361]{SoergelWendt:Perverse}, the restriction of $\R_\calH$ to Tate motives factors as 
$$\R_\calH : \DTM(S)_\C \r \DTM_\calH(S) \stackrel \cong \r \grMod_\C.$$
Since $\R_\calH$ is symmetric monoidal colimit-preserving functor 
it corresponds to 
a morphism
$A \r \C$. By the above remark, it is determined by its restriction along the unit $1 : \C \r A$, which is a simply the identity map $\id_\C$.

On the other hand, by definition, the augmentation map $a : A \r \C$ also has the property that its restriction along the unit map is $\id_\C$.
Therefore, the following two functors are equivalent:
$$\R_\calH \cong (\C \t_{a, A} -) : \DTM(S) \r \grMod_\C.$$
Therefore, there is a functor
\begin{align*}
  \R_{\calH, \red} : \DMr(X) &= \DM(X) \t_{\DTM(S)} \grMod_\C \\&\r \DM_\calH(X) \t_{\DTM_\calH(S)} \grMod_\C = \DM_\calH(X).\qedhere
\end{align*}
\xpf
\prop
\thlabel{comparison.Hodge}
Let $(X,X^+)$ be a scheme with a cellular Whitney--Tate  stratification.
Then there is an equivalence
\begin{center}
  \begin{tikzcd}
  {\DTMr(X,X^+)} \arrow[r,"\R_{\calH, \red}", "\sim"'] & {\DTM_\calH(X,X^+).}
  \end{tikzcd}
\end{center}

\xprop
\pf
Now using that $X$ is Whitney--Tate, we can express $\DTM(X,X^+)$ inductively as a lax limit of a diagram involving the categories $\DTM(X_w)$, as in \thref{DM.lax.limit}.
By \thref{laxlim.tensor}, tensoring with $\grMod_\C$ preserves that lax limit.
Thus, it suffices to prove the claim if $X = \A^n_S \x \GmX S^m$.
By \thref{DTM.flat.cell} $$\DTMr(X) = \Mod_{\Sym(\C\twi{-1}[-1])^{\bigotimes m})}(\grMod_\C).$$
The proof of \thref{DTM.flat.cell} carries over to $\DTM_\calH(X)$, using that $\H^i(\GmX \C^\an, \C) = \C(i)$ for $i=0, 1$ and 0 otherwise.
\xpf

\subsection{Comparison with graded $\ell$-adic sheaves}
Let $S = \Spec \Fq$ be a finite field. We compare reduced motives with the category of graded sheaves introduced very recently by Ho and Li \cite{HoLi:Revisiting}. For simplicity, we restrict our comparison result to the case of schemes, 
referring to \thref{DTM.flat.examples}\refit{stacky.case} for some comments on the case of stacks.

Let $\Shv(X) := \Shv(X,\Qlq)$ be the \ii-category of ind-constructible $\Qlq$-adic sheaves on $X$. Let $\Shv_\mix(X)$ be its full subcategory of ind-mixed complexes, i.e., filtered colimits of mixed complexes as introduced in \cite{BeilinsonBernsteinDeligne:Faisceaux}. 
For example, $\Shv(S)$ is the derived \ii-category of complexes of $\Qlq$-vector spaces equipped with a continuous action of $\Gal(\Fq)$ and the compact objects in $\Shv_{\mix}(S)$ are precisely those perfect complexes where the eigenvalues of $\Frob \in \Gal(\Fq)$ are algebraic numbers whose absolute value is a power of $q^{\frac 12}$. 
The category $\Shv_{\mix}(S)$ decomposes as a coproduct 
of the 
subcategories consisting of those complexes on which $\Frob$ has eigenvalues with absolute value $q^\frac n 2$, for $n \in \Z$. In particular, there is a (colimit preserving, symmetric monoidal) forgetful functor 
$$u : \Shv_{\mix}(S) \r \grMod_\Qlq.$$ 

The category of \emph{graded sheaves} is defined in op.~cit.~as:
$$\Shv_\gr(X) := \Shv_\mix(X) \t_{\Shv_\mix(S)} \grMod_\Qlq,$$
where the tensor product is formed using *-pullback along the structural map $X \r S$, and the above-mentioned forgetful functor. 

\prop
\thlabel{comparison.graded}
Let $X / S$ be a scheme.
Then the $\ell$-adic realization functor $\R_\ell : \DM(X) \r \Shv(X)$ induces a realization functor
$$\R_{\ell, \red} : \DMr(X)_\Qlq \r \Shv_\gr(X).$$
If $X$ is a scheme with a cellular Whitney--Tate  stratification, the restriction of $\R_{\red, \ell}$ is an equivalence 
$$\R_{\ell, \red} : \DTMr(X, X^+)_\Qlq \stackrel \cong \r \Shv_{\Tate, \gr}(X),$$
where the target is the full subcategory of $\Shv_\gr(X)$ consisting of \emph{graded stratified Tate sheaves}, i.e., those graded sheaves $F$ whose restrictions $\iota_w^* F$ to the strata ($\iota_w : X_w \r X$) lie in the presentable subcategory of $\Shv_\gr(X_w)$ generated by the sheaves $\Qlq(n)_{X_w}$ for $n \in \Z$. 
\xprop

\pf Can be shown as in \thref{comparison.Hodge} using the standard calulation for $\H^*(X \x \Spec \Fqq, \Qlq)$ for a single stratum $X=\A^n \x \Gm^m.$
\xpf

\subsection{Comparison with parity sheaves and Achar--Riche's mixed category}
Let $S=\Spec(\C).$ We compare (reduced) stratified Tate motives with \emph{parity sheaves} and Achar--Riche's \emph{mixed category.}

Parity sheaves are certain complexes of sheaves on a complex variety $X^\an$ defined via a condition on the vanishing of the stalk cohomology in even degrees, see \cite{Juteau:ParitySheaves2014}. 
In practice, parity sheaves arise from affine-stratified resolutions of singularities of stratum closures. If $\Char\Lambda=0$, the decomposition theorem implies that parity sheaves are direct sums of intersection cohomology complexes. This is not true if $\Char\Lambda=p$ and parity sheaves often take the role intersection cohomology complexes in modular representation theory, see for example \cite{Soergel:ICandRT, Williamson:AlgebraicRepresentationsConstructible2017}. 

\defi Let $(X,X^+)$ be an (ind-)scheme with an affine Whitney--Tate stratification. Denote by $\Shv(X^{\an})$ the stable \ii-category of sheaves of $\Lambda$-modules on $X^\an.$
The full subcategory of \emph{parity sheaves} $\Par(X^\an,X^+)\subset\Shv(X^{\an})$ consists of all objects $E$ for which the cohomology sheaves on the strata $\mathcal{H}^i(\iota^{*}E)$ are constant, finitely generated, non-zero only in finitely many degrees and zero if $i\not\in2\Z.$
\xdefi
\rema In \cite{Juteau:ParitySheaves2014}, parity sheaves are defined in greater generality.  
Moreover, \cite[Definition 2.4]{Juteau:ParitySheaves2014} actually yields $\Par(X^{\an},X^+)\oplus \Par(X^{\an},X^+)[1].$
\xrema
We compare motives and parity sheaves via the \emph{Betti realisation functor}
$$\BR:\DM(X)\to \Shv(X^{\an})$$
which is compatible with the six functor formalism on both categories. 
\prop \thlabel{DTM.flat.vs.DTM.vs.Par}
Let $(X,X^+)$ be an (ind-)scheme with an affine stratification that admits affine-stratified resolution, see \thref{resolution.condition}.
Then Betti realisation and the reduction functor give equivalences of additive categories
\begin{center}

\begin{tikzcd}[column sep=14pt]
  {\Ho(\Par(X^{\an},X^+))} & {\Ho(\DTM(X,X^+)^{\comp,\wt=0})} \arrow[l, "\BR"',"\sim"] \arrow[r, "\red","\sim"'] & {\Ho(\DTMr(X,X^+)^{\comp,\wt=0})}
  \end{tikzcd}
\end{center}
\xprop 

\pf The statement about the reduction functor is \thref{weight.complex.equivalence}. 
The statement about $\BR$ can with shown with very similar arguments: First, the functor $\BR$ is fully faithful on pointwise pure objects. This can be reduced to the case of single stratum by induction. Second, both the categories of parity sheaves and weight zero reduced stratified Tate motives are generated by pushforwards of constant objects on the resolution of the strata which implies the essential surjectivity.
\xpf
In a nice formalism of mixed sheaves on spaces with an affine stratification, such as the formalism $\DTMr(X, X^+)$ discussed here, the weight complex functor should yield an equivalence between mixed sheaves and the category of chain complexes of weight zero objects. In \cite[Section 7.2]{Achar:KoszulDualitySemisimplicity2011} and \cite[Section 2.2]{Achar:ModularPerverseSheavesII} Achar--Riche take the ingenious approach of simply defining their \emph{mixed category} 
$$\DAR(X^{\an}):=\Ch^\bd(\Ho(\Par(X^{\an},X^+)))$$
via this property. We immediately obtain the following comparison.
\prop \thlabel{DTM.flat.paritysheaves}Under the assumptions of \thref{DTM.flat.vs.DTM.vs.Par} there is an equivalence
\begin{center}
\begin{tikzcd}
{\DTMr(X,X^+)^{\comp}} \arrow[r, "\sim"] & {\DAR(X^\an).}
  \end{tikzcd}
\end{center}
\xprop
\pf Follows from \thref{weight.complex.equivalence} and \thref{DTM.flat.vs.DTM.vs.Par}.  
\xpf
\rema 
In \cite[Section 2.3]{Achar:ModularPerverseSheavesII} Achar--Riche construct pullback and pushforward functors for $\DAR$ in the case of locally closed inclusions of strata as well as affine-stratified proper morphism. The functors are defined via taking compositions and adjoints of functors that preserve parity sheaves and can hence be applied pointwise on $\Ch^\bd(\Ho(\Par(X^\an,X^+))).$
Using that the weight complex functor commutes with weight exact functors \cite{Sosnilo:Theorem} it can be shown that these functors admit a similar description for $\DTMr(X,X^+)$ and are thereby compatible with the comparison in \thref{DTM.flat.paritysheaves}.

\xrema

\bibliographystyle{alphaurl}
\bibliography{bib}

\end{document}